\documentclass{amsart}

\usepackage{amsthm}
\usepackage{amsfonts}
\usepackage{amsmath}
\usepackage{amssymb}
\usepackage{mathrsfs}
\usepackage{fullpage}
\usepackage{stmaryrd}
\usepackage{hyperref}
\usepackage{verbatim}

\theoremstyle{plain}

\theoremstyle{definition}

\theoremstyle{remark}

\newcommand{\Begin}[2]{\begin{#1}\label{#2}}

\newcommand{\bPi}{\mathbf{\Pi}}
\newcommand{\bSigma}{\mathbf{\Sigma}}
\newcommand{\bDelta}{\mathbf{\Delta}}

\newcommand{\bbR}{\mathbb{R}}

\newcommand{\bbO}{\mathbb{O}}

\newcommand{\forces}{\Vdash}
\newcommand{\analytic}{{\bSigma_1^1}}
\newcommand{\lanalytic}{{\Sigma_1^1}}
\newcommand{\coanalytic}{{\bPi_1^1}}

\newcommand{\borel}{{\bDelta_1^1}}

\newcommand{\cantorspace}{{{}^\omega 2}}

\newcommand{\finBinarySequence}{{{}^{<\omega}2}}

\newcommand{\ZFC}{\mathsf{ZFC}}
\newcommand{\ZF}{\mathsf{ZF}}

\newcommand{\reals}{\bbR}

\newcommand{\OD}{\mathrm{OD}}
\newcommand{\HOD}{\mathrm{HOD}}
\newcommand{\AD}{\mathsf{AD}}
\newcommand{\AC}{\mathsf{AC}}
\newcommand{\ON}{\mathrm{ON}}

\newcommand{\degrees}{\mathcal{D}}

\begin{document}

\title{Cardinality of Wellordered Disjoint Unions of Quotients of Smooth Equivalence Relations}

\author{William Chan}
\address{Department of Mathematics, University of North Texas, Denton, TX 76203}
\email{William.Chan@unt.edu}

\author{Stephen Jackson}
\address{Department of Mathematics, University of North Texas, Denton, TX 76203}
\email{Stephen.Jackson@unt.edu}

\begin{abstract}
Assume $\mathsf{ZF + AD^+ + V = L(\mathscr{P}(\reals))}$. Let $\approx$ denote the relation of being in bijection. Let $\kappa \in \ON$ and $\langle E_\alpha : \alpha < \kappa\rangle$ be a sequence of equivalence relations on $\reals$ with all classes countable and for all $\alpha < \kappa$, $\reals \slash E_\alpha \approx \reals$. Then the disjoint union $\bigsqcup_{\alpha < \kappa} \reals \slash E_\alpha$ is in bijection with $\reals \times \kappa$ and $\bigsqcup_{\alpha < \kappa} \reals \slash E_\alpha$ has the J\'onsson property. 

Assume $\mathsf{ZF + AD^+ + V = L(\mathscr{P}(\reals))}$. A set $X \subseteq [\omega_1]^{<\omega_1}$ has a sequence $\langle E_\alpha : \alpha < \omega_1\rangle$ of equivalence relations on $\reals$ such that $\reals \slash E_\alpha \approx \reals$ and $X \approx \bigsqcup_{\alpha < \omega_1} \reals \slash E_\alpha$ if and only if $\reals \sqcup \omega_1$ injects into $X$. 

Assume $\AD$. Suppose $R \subseteq [\omega_1]^\omega \times \reals$ is a relation such that for all $f \in [\omega_1]^\omega$, $R_f = \{x \in \reals : R(f,x)\}$ is nonempty and countable. Then there is an uncountable $X \subseteq \omega_1$ and function $\Phi : [X]^\omega \rightarrow \reals$ which uniformizes $R$ on $[X]^\omega$: that is, for all $f \in [X]^\omega$, $R(f,\Phi(f))$.

Under $\AD$, if $\kappa$ is an ordinal and $\langle E_\alpha : \alpha < \kappa\rangle$ is a sequence of equivalence relations on $\reals$ with all classes countable, then $[\omega_1]^\omega$ does not inject into $\bigsqcup_{\alpha < \kappa} \reals \slash E_\alpha$. 
\end{abstract}

\thanks{March 6, 2019. The first author was supported by NSF grant DMS-1703708. The second author was supported by NSF grant DMS-1800323.}

\maketitle


\section{Introduction}\label{introduction}

The original motivation for this work comes from the study of a simple combinatorial property of sets using only definable methods. The combinatorial property of concern is the J\'onsson property: Let $X$ be any set. For each $n \in \omega$, let $[X]^n_= = \{f \in {}^nX : (\forall i,j \in n)(i \neq j \Rightarrow f(i) \neq f(j))\}$. Let $[X]^{<\omega}_=  = \bigcup_{n \in \omega} [X]^n_=$. A set $X$ has the J\'onsson property if and only if for every function $F : [X]^{<\omega}_= \rightarrow X$, there is some $Y$ such that $Y \approx X$ (where $\approx$ denotes the bijection relation) so that $F[[Y]^{<\omega}_=] \neq X$. That is, $F$ can be made to miss at least one point in $X$ when restricted to the collection of finite unequal tuples of some subset $Y$ of $X$ of the same cardinality as $X$.

Under the axiom of choice, if there is a set with the J\'onsson property, then large cardinal principles such as $0^\sharp$ hold. Using a measurable cardinal, one can construct models of $\ZFC$ in which $2^{\aleph_0}$ is J\'onsson and is not J\'onsson. Hence assuming the consistency of some large cardinals, the J\'onsson property of $2^{\aleph_0}$ is independent of $\ZFC$. Using $\AC$, the sets $\reals$, $\reals \sqcup \omega_1$, $\reals \times \omega_1$, and $\reals \slash E_0$ are all in bijection. ($E_0$ is the equivalence relation defined on $\reals = \cantorspace$ by $x \ E_0 \ y$ if and only if $(\exists m)(\forall n \geq m)(x(n) = y(n))$.)

From a definability perspective, the sets $\reals$, $\reals \times \omega_1$, $\reals \sqcup \omega_1$, and $\reals \slash E_0$ do not have definable bijections without invoking definable wellorderings of the reals which can exist  in canonical inner models like $L$ but in general can not exist if the universe satisfies more regularity properties for sets of reals. For example, there are no injections of $\reals \slash E_0$ into $\reals$ that is induced by a $\borel$ reduction $\Phi : \reals \rightarrow \reals$ of the $=$ relation into the $E_0$ equivalence relation. Such results for the low projective pointclasses can be extended to all sets assuming the axiom of determinacy, $\AD$. Methods that hold in a determinacy setting are often interpreted to be definable methods. Moreover, the extension of $\AD$ called $\AD^+$ captures this definability setting even better since $\AD^+$ implies all sets of reals have a very absolute definition known as the $\infty$-Borel code.

Kleinberg \cite{Infinitary-Combinatorics-and-the-Axiom-of-Determinateness} showed that $\aleph_n$ has the J\'onsson property for all $n \in \omega$ under $\AD$. \cite{Determinacy-and-Jonsson-Cardinals-in-LR} showed that under $\ZF + \AD + V = L(\reals)$, every cardinal below $\Theta$ has the J\'onsson property. (Woodin showed that $\AD^+$ alone can prove this result.)

Holshouser and Jackson began the study of the J\'onsson property for nonwellorderable sets under $\AD$ such as $\reals$. In \cite{Partition-Properties-for-Non-Ordinal-Sets}, they showed that $\reals$ and $\reals \sqcup \omega_1$ have the J\'onsson property. They also showed, using that fact that all $\kappa < \Theta$ have the J\'onsson property, that $\reals \times \kappa$ is J\'onsson. \cite{Definable-Combinatorics-Some-Borel-Equivalence-Relations-arxiv} showed that $\reals \slash E_0$ does not have the J\'onsson property.

Holshouser and Jackson then asked if the J\'onsson property of sets is preserved under various operations. The disjoint union operation will be the main concern of this paper: If $\kappa \in \mathrm{ON}$ and $\langle X_\alpha : \alpha < \kappa\rangle$ is a sequence of sets with the J\'onsson property, then does the disjoint union $\bigsqcup_{\alpha < \kappa} X_\alpha$ have the J\'onsson property? (Here $\bigsqcup$ will always refer to a formal disjoint union in constract to the ordinary union $\bigcup$.) More specifically, does a disjoint union of sets, each in bijection with $\reals$, have the J\'onsson property? The determinacy axioms are particular helpful for studying sets which are surjective images of $\reals$. Hence, a natural question would be if $\langle E_\alpha : \alpha < \kappa\rangle$ is a sequence of equivalence relations on $\reals$ such that for each $\alpha$, $\reals \slash E_\alpha$ is in bijection with $\reals$, then does $\bigsqcup_{\alpha < \kappa} \reals \slash E_\alpha$ have the J\'onsson property? An equivalence relation $E$ on $\reals$ is called smooth if and only if if $\reals \slash E$ is in bijection with $\reals$. (Note that this term is used differently than the ordinary Borel theory which would define $E$ to be smooth if $\reals \slash E$ injects into $\reals$. This will be refered as being weakly smooth.)

$\borel$ equivalence relations with all classes countable are very important objects of study in classical invariant descriptive set theory. One key property that make their study quite robust is the Lusin-Novikov countable section uniformization, which for instance, can prove the Feldman-Moore theorem. \cite{Ordinal-Definability-and-Combinatorics-Equivalence-Relations-arxiv} attempted to study the J\'onsson property for disjoint unions of smooth equivalence relations with all classes countable. It was shown in \cite{Ordinal-Definability-and-Combinatorics-Equivalence-Relations-arxiv} that if $\langle E_\alpha : \alpha < \kappa\rangle$ is a sequence of equivalence relations with all classes countable (not necessarily smooth) and $F : [\bigsqcup_{\alpha < \kappa} \reals \slash E_\alpha]^{<\omega}_= \rightarrow \bigsqcup_{\alpha < \kappa} \reals \slash E_\alpha$, then there is a perfect tree $p$ on $2$ so that 
$$F\left[[\bigsqcup_{\alpha < \kappa} [p] \slash E_\alpha]^{<\omega}_=\right] \neq \bigsqcup_{\alpha < \kappa} \reals \slash E_\alpha.$$
(Here $\reals$ refers to the Cantor space $\cantorspace$.) This ``psuedo-J\'onsson property'' would imply the true J\'onsson property if $\bigsqcup_{\alpha < \kappa} [p] \slash E_\alpha$ is in bijection with $\bigsqcup_{\alpha < \kappa} \reals \slash E_\alpha$. In general for nonsmooth equivalence relations, this can not be true since, for example, $E_0$ is an equivalence relation with all classes countable and $\reals \slash E_0$ is not J\'onsson (\cite{Definable-Combinatorics-Some-Borel-Equivalence-Relations-arxiv}). When each $E_\alpha$ is the identity relation $=$, then one can demonstrate these two sets are in bijection. This (\cite{Ordinal-Definability-and-Combinatorics-Equivalence-Relations-arxiv}) shows that $\reals \times \kappa$ has the J\'onsson property, where $\kappa$ is any ordinal, using only classical descriptive set theoretic methods and does not rely on any combinatorial properties of the ordinal $\kappa$. 

\cite{Ordinal-Definability-and-Combinatorics-Equivalence-Relations-arxiv} asked if $\langle E_\alpha : \alpha < \kappa\rangle$ consists entirely of smooth equivalence relations on $\reals$ with all classes countable and $p$ is any perfect tree on $2$, then is $\bigsqcup_{\alpha < \kappa} \reals \slash E_\alpha$ and $\bigsqcup_{\alpha < \kappa} [p]\slash E_\alpha$ in bijection? Do such disjoint unions have the J\'onsson property? The most natural attempt to show that wellordered disjoint unions of quotients of smooth equivalence relations with all classes countable is J\'onsson would be to show it is, in fact, in bijection with $\reals \times \kappa$, which has already been shown to possess the J\'onsson property.

The computation of the cardinality of wellordered disjoint unions of quotients of smooth equivalence relations on $\reals$ with all classes countable is the main result of the paper. Under $\AD^+$, any equivalence relation on $\reals$ has an $\infty$-Borel code. However, for the purpose of this paper, given a sequence of equivalence relations $\langle E_\alpha : \alpha < \kappa\rangle$ on $\reals$, one will need to uniformly obtain $\infty$-Borel codes for each $E_\alpha$. It is unclear if this is possible under $\AD^+$ alone. For the purpose of obtaining this uniformity of $\infty$-Borel codes, one will need to work with natural models of $\AD^+$, i.e. the axiom system $\mathsf{ZF + AD^+ + V = L(\mathscr{P}(\reals))}$. 

It should be noted that the assumption that each equivalence relation has all classes countable is very important. Analogous to the role of the Lusin-Novikov countable section uniformization in the classical setting, Woodin's countable section uniformization under $\AD^+$ will play a crucial role. 

There are some things that can be said about $\bigsqcup_{\alpha < \kappa} \reals \slash E_\alpha$ when $\langle E_\alpha : \alpha < \kappa\rangle$ is a sequence of smooth equivalence relation (with possibily uncountable classes). It is immediate that $\bigsqcup_{\alpha < \kappa} \reals \slash E_\alpha$ will contain a copy of $\omega_1 \sqcup \reals$. Hence $\omega_1 \sqcup \reals$ is a lower bound on the cardinality of disjoint unions of quotients of smooth equivalence relations. An example of Holshouser and Jackson (Fact \ref{disjoint union smooth equal reals union omega1}) produces a sequence $\langle F_\alpha : \alpha < \omega_1\rangle$ of smooth equivalence relations such that $\bigsqcup_{\alpha < \omega_1} \reals \slash F_\alpha$ is in bijection with $\omega_1 \sqcup \reals$. So this lower bound is obtainable. In fact, in natural models of $\AD^+$, if a set $X \subseteq [\omega_1]^{<\omega_1}$ contains a copy of $\reals \sqcup \omega_1$, then it can be written as an $\omega_1$-length disjoint union of quotients of smooth equivalence relations:
\\*
\\*\textbf{Theorem \ref{characterization disjoint union smooth}} \textit{Assume $\mathsf{ZF + AD^+ + V = L(\mathscr{P}(\reals))}$. Suppose $X \subseteq [\omega_1]^{<\omega_1}$ and $\reals \sqcup \omega_1$ injects into $X$. Then there exists a sequence $\langle E_\alpha : \alpha < \omega_1\rangle$ of smooth equivalence relations on $\reals$ so that $X$ is in bijection with $\bigsqcup_{\alpha < \omega_1} \reals \slash E_\alpha$.}

\textit{Therefore, $X \subseteq [\omega_1]^{<\omega_1}$ has a sequence $\langle E_\alpha : \alpha < \omega_1\rangle$ of smooth equivalence relations such that $X \approx \bigsqcup_{\alpha < \omega_1} \reals \slash E_\alpha$ if and only if $\reals \sqcup \omega_1$ injects into $X$.}
\\*
\\*\indent  $\reals \times \kappa$ is a disjoint union coming from $\langle E_\alpha : \alpha < \kappa\rangle$ where each $E_\alpha$ is the $=$ relation, which is an equivalence relation with all classes countable. However, the proof of Theorem \ref{characterization disjoint union smooth} uses equivalence relations with uncountable classes. Intuitively, it seems that $\reals \sqcup \omega_1$, $[\omega_1]^\omega$, and $[\omega_1]^{<\omega_1}$ should not be obtainable using equivalence relations with countable classes. This motivates the conjecture that if $\langle E_\alpha : \alpha < \kappa\rangle$ is a sequence of smooth equivalence relations with all classes countable then the cardinality of $\bigsqcup_{\alpha < \kappa} \reals \slash E_\alpha$ is $\reals \times \kappa$.

Woodin \cite{The-Cardinals-Below-CountableSubsetOmega} showed that there is elaborate structure of cardinals below $[\omega_1]^{<\omega_1}$. Theorem \ref{characterization disjoint union smooth} shows that every cardinal above $\omega_1 \sqcup \reals$ and below $[\omega_1]^{<\omega_1}$ is a disjoint union of quotients of smooth equivalence relations. A success in Holshouser and Jackson's original goal of establishing the closure of the J\'onsson property under disjoint union would yield the J\'onsson property for many cardinals below $[\omega_1]^{<\omega_1}$. On the other hand, it difficult to see how one could establish the J\'onsson property for every set that appears in this rich cardinal structure using solely the manifestation of these sets as a disjoint union of quotients of smooth equivalence relations.

The main tool for computing the cardinal of wellordered disjoint union of quotients of smooth equivalence relations with all countable classes is the Woodin's perfect set dichotomy which generalizes the Silver's dichotomy for $\coanalytic$ equivalence relations. This perfect set dichotomy states that under $\AD^+$, for any equivalence relation $E$ on $\reals$, either (i) $\reals \slash E$ is wellorderable or (2) $\reals$ injects into $\reals \slash E$. Section \ref{perfect set dichotomy and wellorderable section uniformization} is dedicated to proving this result. A detailed analysis of the proof of this result will be needed for this paper. The proof for case (i) yields a uniform procedure which takes an $\infty$-Borel code for $E$ and gives a wellordering of $\reals \slash E$. Moreover it shows that in this case, $\reals \slash E \subseteq \OD_S$, where $S$ is the  $\infty$-Borel code for $E$. Under $\mathsf{ZF + AD^+ + V = L(\mathscr{P}(\reals))}$, this will give a more general countable section uniformization (Fact \ref{wellorderable section uniformization}). It will be seen that in the proof of case (ii), the injection of $\reals$ will depend on certain parameters. If these parameters could be found uniformly for each equivalence relation from the sequence $\langle E_\alpha : \alpha < \kappa\rangle$, then the proof in case (ii) can uniformly produce injections of $\reals$ into each $\reals \slash E_\alpha$. Together, one would get an injection of $\reals \times \kappa$ into $\bigsqcup_{\alpha < \kappa} \reals \slash E_\alpha$. In general this can not be done; for instance using the example from Fact \ref{disjoint union smooth equal reals union omega1}. However, this can be done when all the equivalence relations are smooth and have all classes countable then one can prove the following:
\\*
\\*\textbf{Theorem \ref{countable smooth R times kappa injection}} \textit{Assume $\mathsf{ZF + AD^+ + V = L(\mathscr{P}(\reals))}$. Let $\kappa \in \ON$ and $\langle E_\alpha : \alpha < \kappa\rangle$ be a sequence of smooth equivalence relations on $\reals$ with all classes countable. Then $\reals \times \kappa$ injects into $\bigsqcup_{\alpha < \kappa} \reals \slash E_\alpha$.}
\\*
\\*\noindent This shows that $\reals \times \kappa$ is a lower bound for the cardinal of $\bigsqcup_{\alpha < \kappa} \reals \slash E_\alpha$. 

Section \ref{upper bound on cardinality} will provide the proof of the relevant half of Hjorth's generalized $E_0$-dichotomy. Again, what is important from this result is the observation that if $\reals \slash E_0$ does not inject into $\reals \slash E$, then there is a wellordered separating family for $E$ defined uniformly from the $\infty$-Borel code for $E$. If $\mathsf{ZF + AD^+ + V = L(\mathscr{P}(\reals))}$ holds, then one has a uniform sequence of $\infty$-Borel codes for the sequence of equivalence relations $\langle E_\alpha : \alpha < \kappa\rangle$, where each $E_\alpha$ is smooth. Using the argument of Hjorth's dichotomy, one obtains uniformly a separating family for each $E_\alpha$. This gives a sequence of injections of each $\reals \slash E_\alpha$ into $\mathscr{P}(\delta)$ where $\delta$ is a possibly very large ordinal. If $\langle E_\alpha : \alpha < \kappa\rangle$ consists entirely of equivalence relations with all classes countable, then the generalized countable section uniformization can be used to uniformly obtain a selector for each $\reals \slash E_\alpha$. This gives the desired injection into $\reals \times \kappa$:
\\*
\\*\textbf{Theorem \ref{upper bound disjoint union smooth}} \textit{ Assume $\mathsf{ZF + AD^+ + V = L(\mathscr{P}(\reals))}$. Let $\kappa$ be an ordinal and $\langle E_\alpha : \alpha < \kappa\rangle$ be a sequence of smooth equivalence relations on $\reals$ with all classes countable. Then there is an injection of $\bigsqcup_{\alpha < \kappa} \reals \slash E_\alpha$ into $\reals \times \kappa$.}
\\*
\\*\textbf{Theorem \ref{cardinality disjoint union smooth}} \textit{ Assume $\mathsf{ZF + AD^+ + V = L(\mathscr{P}(\reals))}$. Let $\kappa \in \ON$ and $\langle E_\alpha : \alpha < \kappa\rangle$ be a sequence of smooth equivalence relations on $\reals$ with all classes countable. Then $\bigsqcup_{\alpha < \kappa} \reals \slash E_\alpha \approx \reals \times \kappa$ and hence $\bigsqcup_{\alpha < \kappa} \reals \slash E_\alpha$ has the J\'onsson property.}
\\*
\\*\indent Many of the results of \cite{Ordinal-Definability-and-Combinatorics-Equivalence-Relations-arxiv} concerning the J\'onsson property of disjoint unions of quotients of equivalence relations on $\reals$ with all classes countable were originally proved under $\AD^+$ using (full) countable section uniformization for relations on $[\omega_1]^\omega \times \reals$. However, most of the results held in merely $\AD$ by using just a form of almost-full uniformization, for example comeager uniformization for relations on $\reals \times \reals$. 

$[\omega_1]^\omega$ is the collection of increasing functions $f : \omega \rightarrow \omega_1$. \cite{Ordinal-Definability-and-Combinatorics-Equivalence-Relations-arxiv} proved under $\AD^+$ that $[\omega_1]^\omega$ does not inject into $\bigsqcup_{\alpha < \kappa} \reals \slash E_\alpha$ when $\langle E_\alpha : \alpha < \kappa\rangle$ is a sequence of equivalence relations on $\reals$ with all classes countable. (Of course, Theorem \ref{cardinality disjoint union smooth} asserts that such a disjoint union is in bijection with $\reals \times \kappa$ under the assumption $\mathsf{ZF + AD^+ + V = L(\mathscr{P}(\reals))}$.) The key ingredient is the ability to uniformize relations $R \subseteq [\omega_1]^\omega \times \reals$ such that for all $f \in [\omega_1]^\omega$, $R_f = \{x \in \reals : R(f,x)\}$ is nonempty and countable. Such a full uniformization is provable under $\AD^+$. A careful inspection of the argument will show that one only needs to uniformize this relation on some $Z \subseteq [\omega_1]^\omega$ such that $Z \approx [\omega_1]^\omega$ to show that no such injection exists. Question 4.21 of \cite{Ordinal-Definability-and-Combinatorics-Equivalence-Relations-arxiv} asks whether such an almost full uniformization is provable in $\AD$. 

It should be noted that if one drops the demand that $R_f$ be countable, then one cannot prove this in general. (See the discussion in Section \ref{almost full countable section uniformization}.) The final section will show that such an almost full countable section uniformization for relations on $[\omega_1]^\omega \times \reals$ is provable in $\AD$:
\\*
\\*\textbf{Theorem \ref{representation of relation on seq}} \textit{$(\ZF + \AD)$ Let $R \subseteq [\omega_1]^\omega \times \reals$ be such that for all $f \in [\omega_1]^\omega$, $R_f$ is nonempty and countable. Then there exists some uncountable $X \subseteq \omega_1$ and function $\Psi$ which uniformizes $R$ on $[X]^\omega$: For $f \in [X]^\omega$, $R(f,\Psi(f))$.}
\\*
\\*\textbf{Corollary \ref{no injection seq into disjoint union}} \textit{$(\ZF + \AD)$ Let $\langle E_\alpha : \alpha < \kappa\rangle$ be a sequence of equivalence relations on $\reals$ with all classes countable, then $[\omega_1]^\omega$ does not inject into $\bigsqcup_{\alpha < \kappa} \reals \slash E_\alpha$.}
\\
\\*\indent These methods also show that for an arbitrary function $\Phi : [\omega_1]^\omega \rightarrow \reals$, one can find some uncountable $X \subseteq \omega_1$ and some reals $\sigma$ and $w$ so that $\Phi(f) \in L[\sigma,w,f]$, for all $f \in [X]^\omega$:
\\*
\\*\textbf{Theorem \ref{constructibility value of function almost everywhere}} \textit{$(\ZF + \AD)$. Let $\Phi : [\omega_1]^\omega \rightarrow \reals$ be a function. Then there is an uncountable $X \subseteq\omega_1$, reals $\sigma,w \in \reals$, and a formula $\phi$ so that for all $f \in [X]^\omega$, $\Phi(f) \in L[\sigma,w,f]$ and for all $z \in \reals$, $z = \Phi(f)$ if and only if $L[\sigma,w,f,z] \models \phi(\sigma,w,f,z)$.}

\section{Basics}

\Begin{definition}{infinity borel codes}
Let $S$ be a set of ordinals and $\varphi$ be a formula in the language of set theory. $(S,\varphi)$ is called an $\infty$-Borel code. For each $n \in \omega$, $\mathfrak{B}^n_{(S,\varphi)} = \{x \in \reals^n : L[S,x] \models \varphi(S,x)\}$ is the subset of $\reals^n$ coded by $(S,\varphi)$. 

A set $A \subseteq \reals^n$ is $\infty$-Borel if and only if $A = \mathfrak{B}^n_{(S,\varphi)}$ for some $\infty$-Borel code $(S,\varphi)$. $(S,\varphi)$ is called an $\infty$-Borel code for $A$.
\end{definition}

\Begin{definition}{AD^+}
(\cite{Axiom-of-Determinacy-Forcing-Axioms} Section 9.1) $\AD^+$ consists of the following statements:

(1) $\mathsf{DC}_\reals$.

(2) Every $A \subseteq \reals$ has an $\infty$-Borel code.

(3) For every $\lambda < \Theta$, $A \subseteq \reals$, and continuous function $\pi: {}^\omega\lambda \rightarrow \reals$, $\pi^{-1}[A]$ is determined.
\end{definition}

Models of the theory $\mathsf{ZF + AD^+ + V = L(\mathscr{P}(\reals))}$ are known as natural models of $\AD^+$. Natural models of $\AD^+$ have several desirable properties. Woodin has shown that these models take one of two forms.

\Begin{fact}{forms of natural models}
(Woodin, \cite{A-Trichotomy-Theorem-in-Natural} Section 3.1) Suppose $\mathsf{ZF + AD^+ + V = L(\mathscr{P}(\reals))}$. If $\AD_\reals$ fails, then there is a set of ordinals $J$ so that $V = L(J,\reals)$. 
\end{fact}

If $V = L(J,\reals)$ for some set of ordinals $J$, then it is clear that every set is ordinal definable from $J$ and a real. That is, every set is ordinal definable from some set of ordinals. This is also true in natural models of $\AD^+$ in which $\AD_\reals$ holds: 

\Begin{fact}{natual models every set od in countable set}
(Woodin, \cite{A-Trichotomy-Theorem-in-Natural} Theorem 3.3) Assume $\mathsf{ZF + AD^+ + AD_\reals + V = L(\mathscr{P}(\reals))}$. Then every set is $\OD$ from some element of $\bigcup_{\lambda < \Theta} \mathscr{P}_{\omega_1}(\lambda)$. 
\end{fact}

\Begin{fact}{location od infinity borel code}
(Woodin, \cite{A-Trichotomy-Theorem-in-Natural} Theorem 3.4) Assume $\mathsf{ZF + AD^+ + V = L(\mathscr{P}(\reals))}$. Let $S$ be a set of ordinals. If $A \subseteq \reals$ is $\OD_S$, then $A$ has an $\OD_S$ $\infty$-Borel code.
\end{fact}

\Begin{definition}{martin measure}
If $x,y \in \reals$, $x \leq_T y$ indicates $x$ is Turing reducible to $y$. $x \equiv_T y$ denotes $x \leq_T y$ and $y \leq_T x$. If $x \in \reals$, then $[x]_T$ denotes the equivalence class of $x$ under $\equiv_T$. Let $\degrees$ denote the collection of $\equiv_T$ equivalence classes.

For $X,Y \in \degrees$, define $X \leq_T Y$ if and only if for all $x \in X$ and $y \in Y$, $x \leq_T y$. $U \subseteq \degrees$ contains a Turing cone with base $X \in \degrees$ if and only if for all $Y \in \degrees$, $X \leq_T Y$ implies that $Y \in U$. 

Let $\mu \subseteq \mathscr{P}(\degrees)$ be defined by $U \in \mu$ if and only if $\mu$ contains a Turing cone. Under $\AD$, $\mu$ is a countably complete ultrafilter on $\degrees$. $\mu$ is called the Martin's measure.
\end{definition}

\Begin{fact}{wellfoundedness martin ultraproduct}
(Woodin, \cite{A-Trichotomy-Theorem-in-Natural} Section 2.2) $\ZF + \AD^+$ proves $\prod_{X \in \degrees} \ON \slash \mu$ is wellfounded.
\end{fact}

\Begin{definition}{vopenka forcing}
For $n \in \omega$ and a set of ordinals $S$, let ${}_n\bbO_S$ denote the collection of nonempty $\OD_S$ subsets of $\reals^n$. (${}_1\bbO_S$ will be denoted by $\bbO_S$.) 

For $p,q \in {}_n\bbO_S$, let $p \leq_{{}_n\bbO_S} q$ if and only if $p \subseteq q$. Let $1_{{}_n\bbO_S} = \reals^n$. $({}_n\bbO_S,\leq_{{}_n\bbO_S}, 1_{{}_n\bbO_S})$ is the $n$-dimensional $S$-Vop\v{e}nka's forcing.

By using an $S$-definable bijection of the collection of $\OD_S$ sets with $\ON$, ${}_n\bbO_S$ can be considered as a set of ordinals. In this way, the forcing ${}_n\bbO_S$ is a forcing belonging to $\HOD_S$. 

For each $m \in \omega$, let $b_m \in \bbO_S$ be defined as $\{x \in \reals : m \in x\}$. Let $\tau = \{(\check m, b_m) : m \in \omega\}$. $\tau$ is an $\bbO_S$-name for a real. (Similar definition exists for all ${}_n\bbO_S$.)
\end{definition}

\Begin{fact}{vopenka theorem}
(Vop\v{e}nka's theorem) Let $M$ be a transitive inner model of $\ZF$. Let $S \in M$ be a set of ordinals. 

For all $x \in \reals^M$, there is an $\bbO_S^M$-generic filter over $\HOD^M_S$, $G_x \in M$, so that $\tau[G_x] = x$. 

Suppose $K$ is an $\OD_S^M$-set of ordinals and $\varphi$ is a formula. Let $N$ be some transitive inner model with $\HOD_S^M \subseteq N$. Suppose $p = \{x \in \reals : L[K,x] \models \varphi(K,x)\}$ is a condition of $\bbO_S^M$ (i.e. is nonempty). Then $N \models p \forces_{\bbO_S^M} L[\check K,\tau] \models \varphi(\check K,\tau)$. 
\end{fact}

\begin{proof}
A proof of this can be found among \cite{Set-Theory} Theorem 15.46, \cite{Dichotomy-for-Definable-Universe} Theorem 2.4, or \cite{L(R)-With-Determinacy-Satisfies-Suslin-Hypothesis} Fact 2.7.
\end{proof}

\Begin{fact}{product vopenka generic}
Let $M$ be an inner model of $\ZF$. Let $S \in M$ be a set of ordinals. Let $N$ be an inner model of $\ZF$ such that $N \supseteq \HOD_S^M$. Let $n \geq 1$ be a natural number. Suppose $(g_0, ..., g_{n - 1})$ is an ${}_n\bbO_S^M$-generic real over $N$. Then each $g_0$, ..., $g_{n-1}$ is $\bbO_S^M$-generic over $N$.
\end{fact}

\begin{proof}
This is straightforward. See \cite{L(R)-With-Determinacy-Satisfies-Suslin-Hypothesis} Fact 2.8 for details.
\end{proof}

\Begin{definition}{cardinality and bijection}
Suppose $X$ and $Y$ are sets. $X \approx Y$ indicates that there is a bijection between $X$ to $Y$. 
\end{definition}

\Begin{definition}{smooth equivalence relation}
Let $E$ be an equivalence relation on $\reals$. $E$ is smooth if and only if $\reals \slash E \approx \reals$. $E$ is weakly smooth if and only if $\reals \slash E$ injects into $\reals$. 

Under $\AD$ by the perfect set property, $E$ is weakly smooth if and only if $\reals \slash E$ is either countable or in bijection with $\reals$. 

(From the theory of Borel equivalence relations, ``smooth'' would usually refer to what is called weakly smooth above. In this article, one will reserve the term ``smooth'' for equivalence relations on $\reals$ whose quotients are in bijection with $\reals$.)
\end{definition}

\section{Perfect Set Dichotomy and Wellorderable Section Uniformization}\label{perfect set dichotomy and wellorderable section uniformization}

The Silver's dichotomy states the every $\coanalytic$ equivalence relation $E$ on $\reals$ has countably many equivalence classes ($\reals \slash E$ is hence wellorderable) or $E$ has a perfect set of pairwise $E$-inequivalent elements ($\reals$ injects into $\reals \slash E$). Woodin's perfect set dichotomy states: under $\AD^+$, for every equivalence relation $E$ on $\reals$, either $\reals \slash E$ is wellorderable or $E$ has a perfect set of $E$-inequivalent elements. As a consequence, every set which is a surjective image of $\reals$, either the set contains a copy of $\reals$ or is wellorderable. In natural models of $\AD^+$, \cite{A-Trichotomy-Theorem-in-Natural} showed that every set either has a copy of $\reals$ or is wellorderable. Moreover, a consequence of the proof shows roughly that every wellorderable $\OD$ set contains only $\OD$ elements. This immediately yields wellorderable section uniformization for rather general relations (on sets) with each section wellorderable. This generalizes Woodin's countable section uniformization for relations on $\reals \times \reals$. 

This section will provide a proof of Woodin's perfect set dichotomy and the wellorderable section uniformization. An observation about the uniformity of the proof of the perfect set dichotomy will be necessary for studying disjoint unions of quotients of smooth equivalence relations with all classes countable. Later countable section uniformization (for relations on $\mathscr{P}(\delta) \times \reals$ where $\delta \in \ON$) will also be needed. All results found in this section are due to Woodin or the authors of \cite{A-Trichotomy-Theorem-in-Natural}.

\Begin{definition}{component of equivalence relation}
Let $E$ be an equivalence relation on $\reals$. An $E$-component is a nonempty set $A$ so that for all $x,y \in A$, $x \ E \ y$. (An $E$-component is just a nonempty subset of an $E$-class.)
\end{definition}

\Begin{theorem}{perfect set dichotomy}
(Woodin) Assume $\ZF + \AD^+$. Let $E$ be an equivalence relation on $\reals$. Then either

(i) $\reals \slash E$ is wellorderable.

(ii) $\reals$ injects into $\reals \slash E$. 
\end{theorem}

\begin{proof}
Silver \cite{Counting-the-Number-of-Equivalence-Classes} proved the $\coanalytic$ version of this result. Harrington \cite{A-Powerless-Proof-of-a-Theorem-of-Silver} produced a proof of this result using the Gandy-Harrington forcing of nonempty $\lanalytic$ subsets of $\reals$. This proof will replace Gandy-Harrington forcing with the Vop\v{e}nka forcing of nonempty $\OD$ subsets of $\reals$. Arguments from \cite{Borel-Equivalence-Relations} Chapter 10 and \cite{Dichotomy-for-Definable-Universe} will be adapted.

Using $\AD^+$, let $(S,\varphi)$ be an $\infty$-Borel code for $E$. In all models considered in this proof, $E$ will always be understood to be the set defined by $(S,\varphi)$.

Note that if $x \equiv_T y$, then $L[S,x] = L[S,y]$, $\HOD^{L[S,x]}_S = \HOD^{L[S,y]}_S$, and the canonical global wellordering of $\HOD_S^{L[S,x]}$ and $\HOD^{L[S,y]}_S$ are the same. Therefore, for each $X \in \degrees$, let $L[S,X] = L[S,x]$ and $\HOD_S^{L[S,X]} = \HOD_S^{L[S,x]}$ for any $x \in X$.

(Case I) For all $X \in \degrees$, for all $a \in \reals^{L[S,X]}$, there is an $\OD^{L[S,X]}_S$ $E$-component $A \in \bbO_S^{L[S,X]}$ containing $a$.

For each $F \in \prod_{X \in \degrees} \omega_1 \slash \mu$, define $A_F$ as follows: Let $f : \degrees \rightarrow \omega_1$ be such that $f \in F$, i.e. is a representative for $F$ under the relation $\sim$ of $\mu$-almost equality. For $a \in \reals$, $a \in A_F$ if and only if on a Turing cone of $X \in \degrees$, $a$ belongs to the $f(X)^\text{th}$ $E$-component in $\bbO_S^{L[S,X]}$ according to the canonical global wellordering of $\HOD^{L[S,X]}_S$. (This $f(X)^\text{th}$ set is said to be $\emptyset$ if there is no $f(X)^\text{th}$ $E$-component in $\bbO_S^{L[S,X]}$.) $A_F$ is well defined in the sense that it is independent of the chosen representative.

$A_F$ is an $E$-component: Suppose $a,b \in A_F$ and $\neg(a \ E \ b)$. Pick $f \in F$. Since $a,b \in A_F$, there is some $Z \geq_T [a \oplus b]_T$ so that for all $X \geq_T Z$, $a$ and $b$ belong to the $f(X)^\text{th}$ $E$-component in $\bbO_S^{L[S,X]}$. This would imply that $a \ E \ b$ holds in $L[S,X]$. Since $E$ is defined by the $\infty$-Borel code $(S,\varphi)$, $L[S,X] \models L[S,a,b] \models \varphi(S,a,b)$. Then $V \models L[S,a,b] \models \varphi(S,a,b)$. Hence $V \models a \ E \ b$. Contradiction.

For all $a \in \reals$, $a$ belongs to some $A_F$: Let $f : \degrees \rightarrow \omega_1$ be defined as follows: For all $X \geq_T [a]_T$, let $f(X)$ be the least $\alpha$ so that $a$ belongs to the $\alpha^\text{th}$ $E$-component of $\bbO_S^{L[S,X]}$. Such an $\alpha$ exists by the Case I assumption. If $F = [f]_\sim$, then $a \in A_F$. 

By Fact \ref{wellfoundedness martin ultraproduct}, $\prod_{X \in \degrees} \omega_1 \slash \mu$ is wellfounded. Hence $\langle A_F : F \in \prod_{X \in \degrees} \omega_1 \slash \mu\rangle$ is a wellordered sequence of $E$-components so that every $a \in \reals$ belongs to some $A_F$. 

Let $B_F$ be the $E$-closure of $A_F$. $\langle B_F : F \in \prod_{X \in \degrees}\omega_1 \slash \mu\rangle$ is a surjection of a wellordered set onto the collection of $E$-classes. By removing duplicates by canonically choosing the least index for each $E$-class, one obtains a bijection $\langle C_\alpha : \alpha < \delta\rangle$ of some $\delta \in \ON$ onto the collection of $E$-classes. Note for later purposes that $\langle C_\alpha : \alpha < \delta\rangle$ is obtained uniformly from the $\infty$-Borel code $(S,\varphi)$. Moreover, each $C_\alpha$ is $\OD_S$.

(Case II) There exists an $X \in \degrees$ and $a \in \reals^{L[S,X]}$ which does not belong to any $\OD_S^{L[S,X]}$ $E$-component of $\bbO_S^{L[S,X]}$.

Let $u$ be the set of $x \in \reals^{L[S,X]}$ that does not belong to any $E$-component in $\bbO_S^{L[S,X]}$. The set $u$ is nonempty by the case II assumption and is $\OD_S^{L[S,X]}$. Hence $u \in \bbO_S^{L[S,X]}$. 

Claim 1: $\HOD_S^{L[S,X]} \models (u,u) \forces_{\bbO_S^{L[S,X]}\times\bbO_S^{L[S,X]}} \neg(\tau_L \ E \ \tau_R)$, where $\tau_L$ and $\tau_R$ are the canonical $\bbO_S^{L[S,X]}\times\bbO_S^{L[S,X]}$-names for the evaluation of the $\bbO_S^{L[S,X]}$-name $\tau$ according the left and right $\bbO_S^{L[S,X]}$-generic coming from an $\bbO_S^{L[S,X]}\times\bbO_S^{L[S,X]}$-generic.

To see this, suppose not. Then there is some $(v,w) \leq_{\bbO_S^{L[S,X]} \times \bbO_S^{L[S,X]}} (u,u)$ so that $\HOD_S^{L[S,X]} \models (v,w) \forces_{\bbO_S^{L[S,X]}\times\bbO_S^{L[S,X]}} \tau_L \ E \ \tau_R$. 

Subclaim 1.1: If $G_0$ and $G_1$ are $\bbO_S^{L[S,X]}$-generic over $\HOD_S^{L[S,X]}$ containing $v$ (but not necessarily mutually generic), then $\tau[G_0] \ E \ \tau[G_1]$. 

To prove Subclaim 1.1: Since $\AD$ implies $\omega_1^V$ is inaccessible in every inner model of $\mathsf{ZFC}$, $\bbO_S^{L[S,X]}$ and its power set in $\HOD_S^{L[S,X]}$ is countable in $V$. Hence, for any $G_0$ and $G_1$ which are $\bbO_S^{L[S,X]}$-generic filters over $\HOD_S^{L[S,X]}$ containing the condition $v$, there exists an $H$ which is $\bbO_S^{L[S,X]}$-generic over $\HOD_S^{L[S,X]}[G_0]$ and $\HOD_{S}^{L[S,X]}[G_1]$ containing the condition $w$. Then by the forcing theorem, $\HOD_S^{L[S,X]}[G_0][H] \models \tau[G_0] \ E \ \tau[H]$ and $\HOD_S^{L[S,X]}[G_1][H] \models \tau[G_1] \ E \ \tau[H]$. Since $E$ is defined by the $\infty$-Borel code $(S,\varphi)$, this means 
$$\HOD_S^{L[S,X]}[G_0][H] \models L[S,\tau[G_0],\tau[H]] \models \varphi(S,\tau[G_0],\tau[H])$$ 
$$\HOD_S^{L[S,X]}[G_1][H] \models L[S,\tau[G_1],\tau[H]] \models \varphi(S,\tau[G_1],\tau[H]).$$
But then
$$V \models L[S,\tau[G_0],\tau[H]] \models \varphi(S,\tau[G_0],\tau[H])$$ 
$$V \models L[S,\tau[G_1],\tau[H]] \models \varphi(S,\tau[G_1],\tau[H])$$ 
This means $V \models \tau[G_0] \ E \ \tau[H]$ and $V \models \tau[G_1] \ E \ \tau[H]$. Therefore, $V \models \tau[G_0] \ E \ \tau[G_1]$.

Note that there exists $a,b \in v$ so that $\neg(a \ E \ b)$ since otherwise $v$ would be an $E$-component in $\bbO_S^{L[S,X]}$. Since $v \subseteq u$, this contradicts the definition of $u$. Let $p = (v \times v) \setminus E$. $p$ is a nonempty $\OD_S^{L[S,X]}$ subset of $\reals^2$. Hence $p \in {}_2\bbO_S^{L[S,X]}$. 

Let $\tau^2$ denote the canonical name for the element of $\reals^2$ added by ${}_2\bbO_S^{L[S,X]}$. Let $\tau^2_0$ and $\tau^2_1$ be the canonical name for the first and second coordinate of $\tau^2$, respectively. Let $G$ be ${}_2\bbO_S^{L[S,X]}$-generic over $\HOD_S^{L[S,X]}$ containing $p$. Since $p \leq_{{}_2\bbO_S^{L[S,X]}} q$, $q \in G$. Since $E$ has $(S,\varphi)$ as its $\infty$-Borel code, the condition $q = \{(x,y) \in (\reals^2)^{L[S,X]} : \neg(x \ E \ y)\}$ can be expressed in the form for which the last statement of Fact \ref{vopenka theorem} applies. Hence $\neg(\tau^2_0[G] \ E \ \tau^2_1[G])$. However Fact \ref{product vopenka generic} states that $\tau^2_0[G]$ and $\tau^2_1[G]$ are the canonical reals added by some $\bbO_S^{L[S,X]}$-generic filter over $\HOD_S^{L[S,X]}$. Subclaim 1.1 implies $\tau^2_0[G] \ E \ \tau^2_1[G]$. Contradiction. This proves Claim 1.

Again by $\AD$, $\bbO_S^{L[S,X]} \times \bbO_S^{L[S,X]}$ and its power set are countable in $V$. Fix an enumeration $(D_n : n \in \omega)$ of all the dense open subsets of $\bbO_S^{L[S,X]}\times\bbO_S^{L[S,X]}$ that belong to $\HOD_S^{L[S,X]}$. By always using the canonical wellordering of $\HOD_S^{L[S,X]}$ to make selections, the routine argument canonically produces a perfect tree so that the $\bbO_S^{L[S,X]}$-generics over $\HOD_S^{L[S,X]}$ containing the condition $u$ along any two different paths are mutually $\bbO_S^{L[S,X]}$-generic over $\HOD_S^{L[S,X]}$. Using Claim 1, the collection of reals added by generic filters along paths of the perfect tree forms a perfect set of pairwise $E$-inequivalent reals.

For later purpose of this paper, note that once one chooses an $X \in \degrees$ witnessing Case II and the enumeration $(D_n : n \in \omega)$ of the dense open subsets of $\bbO_S^{L[S,X]} \times \bbO_S^{L[S,X]}$ in $\HOD_S^{L[S,X]}$, the embedding of $\reals$ into $\reals \slash E$ is given by the explicit procedure above.
\end{proof}

\Begin{fact}{wo od set inside hod}
Assume $\ZF + \AD^+$. If $E$ is an equivalence relation on $\reals$ with $\infty$-Borel code $(S,\varphi)$ and $\reals$ does not inject into the quotient $\reals \slash E$, then $\reals \slash E \subseteq \OD_S$. In particular, if $A \subseteq \reals$ is a countable set of reals with $\infty$-Borel code $(S,\varphi)$, then $A \subseteq \HOD_S$. 

Assume $\mathsf{ZF + AD^+ + V = L(\mathscr{P}(\reals))}$. Suppose $S$ is a set of ordinals. If $E$ is an $\OD_S$ equivalence relation on $\reals$ and $\reals$ does not inject into $\reals \slash E$, then $\reals \slash E \subseteq \OD_S$.  In particular, $A \subseteq \reals$ is a countable $\OD_S$ set of reals, then $A \subseteq \HOD_S$. 

Assume $\mathsf{ZF + \AD^+ + V = L(\mathscr{P}(\reals))}$. If $A$ is any $\OD_S$ set such that $\reals$ does not inject into $A$, then $A \subseteq \OD_S$ and hence wellorderable. 
\end{fact}

\begin{proof}
First work in $\AD^+$, the first statement comes from the observation at the end of the Case I argument in Theorem \ref{perfect set dichotomy} that the sequence $\langle C_\alpha : \alpha < \delta\rangle$ is $\OD_S$, produced uniformly from $(S,\varphi)$, and each $C_\alpha \in \OD_S$. 

If $A \subseteq \reals$ is countable with $\infty$-Borel code $(S,\varphi)$, define the equivalence relation $E$ on $\reals$ by 
$$x \ E \ y \Leftrightarrow (x = y) \vee (x,y \notin A).$$
$E$ has an $\infty$-Borel code which is $\OD_S$. Then apply the first result to $E$.

Now work in $\mathsf{ZF + AD^+ + V = L(\mathscr{P}(\reals))}$. By Fact \ref{location od infinity borel code}, every $\OD_S$ set of reals has an $\infty$-Borel code which is $\OD_S$. The result then follows from apply the earlier statements. 

The final statement is not needed in this paper. However, the idea is that in natural models of $\AD^+$, one can break an arbitrary set $A$ into a uniform sequence of subsets which are surjective images of $\reals$. By the assumption that $A$ does not contain a copy of $\reals$, each of the pieces do not contain a copy of $\reals$. Then Theorem \ref{perfect set dichotomy} uniformly gives a wellordering of each piece. These wellorderings are then coherently patched together into a wellordering of the original set $A$. Recall that by Fact \ref{forms of natural models}, natural models of $\AD^+$ either take the form $L(J,\reals)$ for some set of ordinals $J$ or satisfy $\AD_\reals$. This patching for $L(J,\reals)$ is relatively straightforward. In the $\AD_\reals$ case, it is more challenging and uses the unique supercompactness measure on $\mathscr{P}_{\omega_1}(\lambda)$ for each $\lambda < \Theta$. See \cite{A-Trichotomy-Theorem-in-Natural} or \cite{L(R)-With-Determinacy-Satisfies-Suslin-Hypothesis} for the details.
\end{proof}

\Begin{fact}{wellorderable section uniformization}
Assume $\mathsf{ZF + AD^+ + V = L(\mathscr{P}(\reals))}$. Let $\delta$ be an ordinal and $X$ be a set. Suppose $R \subseteq \mathscr{P}(\delta) \times X$ is a relation so that for each $N \in \mathscr{P}(\delta)$, the section $R_N = \{x \in X : (N,x) \in R\}$ is wellorderable, then $R$ has a uniformization.

In particular, if $R \subseteq \mathscr{P}(\delta) \times \reals$ is a relation so that each $R_N$ is countable, then $R$ has a uniformization.
\end{fact}

\begin{proof}
Using Fact \ref{natual models every set od in countable set} and the remarks preceding this fact, $R$ is ordinal definable from some set of ordinals $S$. So each $R_N$ is ordinal definable from the set of ordinals $\langle S,N\rangle$, where $\langle \cdot ,\cdot\rangle$ refers to some fixed way of coding two sets of ordinals ino a single set of ordinals. Then Fact \ref{wo od set inside hod} implies that $R_N \subseteq \OD_{\langle S, N\rangle}$. The canonical wellordering of $\OD_{\langle S,N\rangle}$ gives a canonical wellordering of $R_N$. 

Although the hypothesis states that each $R_N$ is wellorderable without, a priori, a uniform wellordering, one in fact does have a uniform wellordering of each section. The function that selects the least element of $R_N$ using this canonical uniform wellordering of all sections of $R$ is the desired uniformization function.
\end{proof}

\section{Lower Bound on Cardinality}\label{lower bound on cardinality}

The following section gives a lower bound on the cardinality of disjoint unions of smooth equivalence relations on $\reals$. (Later it will be shown that this fact characterizes those subsets of $[\omega_1]^{<\omega_1}$ which are in bijection with a disjoint union of quotients of smooth equivalence relations.)

\Begin{fact}{lower bound}
$(\ZF)$. Let $\kappa$ be an ordinal. Let $\langle G_\alpha : \alpha < \kappa\rangle$ be a sequence of smooth equivalence relations. Then $\reals \sqcup \kappa$ injects into $\bigsqcup_{\alpha < \kappa} \reals \slash G_\alpha$. 
\end{fact}

\begin{proof}
Let $\bar{0} \in \reals$ denote the constant $0$ function. Let $\Phi : \reals \rightarrow \reals \slash G_0$ so that the image of $\Phi$ does not include $[\bar{0}]_{G_0}$. Let $\Psi : \reals \sqcup \kappa \rightarrow \bigsqcup_{\alpha < \kappa} \reals \slash G_\alpha$ be defined by $\Psi(r) = \Phi(r)$ if $r \in \reals$ and $\Psi(\alpha) = [\bar{0}]_{G_\alpha}$ if $\alpha \in \kappa$. $\Psi$ is an injection.
\end{proof}

This lower bound is optimal by the following example.

\Begin{fact}{disjoint union smooth equal reals union omega1}
$(\ZF)$ There is a sequence $\langle F_\alpha : \alpha < \omega_1\rangle$ of smooth equivalence relation such that $\bigsqcup_{\alpha < \omega_1} \reals \slash F_\alpha \approx \reals \sqcup \omega_1$. 

$(\ZF + \AD)$ There is a sequence $\langle F_\alpha : \alpha < \omega_1\rangle$ of smooth equivalence relation such that $\bigsqcup_{\alpha < \omega_1} \reals \slash F_\alpha$ is not in bijection with $\reals \times \omega_1$. 
\end{fact}

\begin{proof}
Let $\mathrm{WO}$ denote the set of reals coding wellorderings. For each $\alpha < \omega_1$, let $\mathrm{WO}_\alpha$ denote the set of reals coding wellorderings of order type $\alpha$. Define $F_\alpha$ by
$$x \ F_\alpha \ y \Leftrightarrow (x \notin \mathrm{WO}_\alpha \wedge y \notin \mathrm{WO}_\alpha) \vee (x = y).$$
All elements of $\mathrm{WO}_\alpha$ form singleton $F_\alpha$-classes, and there is a single uncountable $\OD$ equivalence class consisting of $\reals \setminus \mathrm{WO}_\alpha$. 

Therefore, $\bigsqcup_{\alpha < \omega_1} \reals \slash F_\alpha \approx (\bigcup_{\alpha < \omega_1} \mathrm{WO}_\alpha) \sqcup \omega_1 = \mathrm{WO} \sqcup \omega_1 \approx \reals \sqcup \omega_1$, where the copy of $\omega_1$ comes from the large equivalence class for each $\alpha < \omega_1$. 

The second statement follows from the fact that under $\mathsf{AD}$, $\reals \sqcup \omega_1$ is not in bijection with $\reals \times \omega_1$. 
\end{proof}

Using countable section uniformization for relations on $[\omega_1]^\omega \times \reals$ given by Fact \ref{wellorderable section uniformization}, \cite{Ordinal-Definability-and-Combinatorics-Equivalence-Relations-arxiv} Fact 4.20 shows that $[\omega_1]^\omega$ can not inject into a disjoint union of quotients of smooth equivalence relations with all classes countable assuming $\ZF + \AD^+$. Section \ref{almost full countable section uniformization} will show under $\AD$ alone that $[\omega_1]^\omega$ can not inject into a disjoint union of quotients of equivalence relations with all section countable by proving an almost full countable section uniformization for relations on $[\omega_1]^\omega \times \reals$.

Countable section uniformization seems to be a powerful tool that allows disjoint unions of quotients of smooth equivalence relations with all classes countable to be studied more easily. Each $F_\alpha$ from the sequence $\langle F_\alpha : \alpha < \omega_1\rangle$ from Fact \ref{disjoint union smooth equal reals union omega1} has only one uncountable class. Its disjoint union $\bigsqcup_{\alpha < \omega_1} \reals \slash F_\alpha$ is in bijection with $\reals \sqcup \omega_1$. A natural question asked in \cite{Ordinal-Definability-and-Combinatorics-Equivalence-Relations-arxiv} was whether it is necessary to use equivalence relations with uncountable classes to produce a disjoint union which is in bijection with $\reals \sqcup \omega_1$. It is also natural to ask if it is possible to determine the cardinality of a disjoint unions of quotients of smooth equivalence relation with all classes countable. 

It will be shown later that many subsets of $[\omega_1]^{<\omega_1}$ are disjoint unions of quotients of smooth equivalence relations on $\reals$. In particularly, $[\omega_1]^{<\omega_1}$ itself is an $\omega_1$-length disjoint union of quotients of smooth equivalence relations. Hence having all classes countable is necessary in the results of this paper.

First, it is elucidating to see why a natural attempt to use the argument in Theorem \ref{perfect set dichotomy} Case II is unable to produce a uniform sequence of embeddings of $\reals$ into $\reals \slash F_\alpha$, where $\langle F_\alpha : \alpha < \omega_1\rangle$ is the sequence of equivalence relations from Fact \ref{disjoint union smooth equal reals union omega1}.

Suppose $\langle E_\alpha : \alpha < \kappa\rangle$ is a wellordered sequence of equivalence relations. The most natural presentation of such a sequence is as a relation $R \subseteq \kappa \times \reals \times \reals$ defined by $(\alpha,x,y) \Leftrightarrow x \ E_\alpha \ y$. Under $\AD^+$, each $E_\alpha \subseteq \reals\times\reals$ has an $\infty$-Borel code. For the results of this paper, one will need to uniformly obtain an $\infty$-Borel code for each $E_\alpha$. It is unclear this can be done for any wellordered sequence of equivalence relations under just $\AD^+$. However, in models of $\mathsf{ZF + AD^+ + V = L(\mathscr{P}(\reals))}$, this will be possible. This motivates the following definition.

\Begin{definition}{infinity borel code for relation on ordinals}
Let $\kappa \in \ON$ and $R \subseteq \kappa \times \reals$. An $\infty$-Borel code for $R$ is a pair $(S,\varphi)$ such that $(\alpha,r) \in R$ if and only if $L[S,r] \models \varphi(S,\alpha,r)$. 

A sequence $\langle E_\alpha : \alpha < \kappa\rangle$ of equivalence relations on $\reals$ has an $\infty$-Borel code if and only if the relation $R(\alpha,x,y) \Leftrightarrow x \ E_\alpha \ y$ has an $\infty$-Borel code.

Note that if $(S,\varphi)$ is an $\infty$-Borel code for $\langle E_\alpha : \alpha < \kappa\rangle$, then uniformly from $(S,\varphi)$, one can find formulas $\varphi_\alpha$ so that $(\langle S,\alpha\rangle,\varphi_\alpha)$ is an $\infty$-Borel code for $E_\alpha$ (in the ordinary sense). (Here $\langle S,\alpha\rangle$ is some fixed coding of sets of ordinals so that $S$ and $\alpha$ can be recovered.)
\end{definition}

\Begin{fact}{natural model AD relation ordinals have infinity borel codes}
Under $\mathsf{ZF + AD^+ + V = L(\mathscr{P}(\reals))}$, for every $\kappa \in \ON$ and every relation $R \subseteq \kappa \times \reals$, $R$ has an $\infty$-Borel code.
\end{fact}

\begin{proof}
By Fact \ref{forms of natural models} and Fact \ref{natual models every set od in countable set}, every set is ordinal definable from some set of ordinals. Let $S$ be a set of ordinals so that $R$ is $\OD_S$. Let $R_\alpha = \{x : (\alpha,x) \in R\}$. Each $R_\alpha$ is $\OD_S$. By Fact \ref{location od infinity borel code}, each $R_\alpha$ has an $\infty$-Borel code in $\HOD_S$. Let $(S_\alpha,\varphi_\alpha)$ be the least $\infty$-Borel code for $R_\alpha$ according to the canonical wellordering of $\HOD_S$. Let $U = \{(\alpha,\beta) : \beta \in S_\alpha\}$. Then there is some $\varphi$ so that $(U,\varphi)$ is an $\infty$-Borel code for $R$ in the sense of Definition \ref{infinity borel code for relation on ordinals}.
\end{proof}

Now suppose that $\langle E_\alpha : \alpha < \kappa\rangle$ is a sequence of smooth equivalence relations. Let $(S,\varphi)$ be an $\infty$-Borel code for this sequence. Hence uniformly, $E_\alpha$ has some $\infty$-Borel code $(\langle S,\alpha\rangle, \varphi_\alpha)$. 

By the remark at the end of the proof of Theorem \ref{perfect set dichotomy}, an embedding of $\reals$ into $\reals \slash E_\alpha$ can be produced uniformly from a choice of $X \in \degrees$ so that the Case II assumptions holds and a fixed enumeration $(D_n : n \in \omega)$ of the dense open subsets of $\bbO_{\langle S,\alpha\rangle}^{L[\langle S,\alpha\rangle,X]} \times\bbO_{\langle S,\alpha\rangle}^{L[\langle S,\alpha\rangle,X]}$. Since for any $\alpha$, $L[\langle S, \alpha\rangle,X] = L[S,X]$ and $\bbO_{\langle S,\alpha\rangle}^{L[\langle S,\alpha\rangle,X]} = \bbO_S^{L[S,X]}$, one can drop the $\alpha$. Thus one can uniformly find a sequence of injections of $\reals$ into $\reals \slash E_\alpha$ if one could find a single $X \in \degrees$ that witnesses the Case II assumption for all equivalence relations $E_\alpha$. With such a sequence, one could then inject $\reals \times \kappa$ into $\bigsqcup_{\alpha < \kappa} \reals \slash E_\alpha$. 

Now consider the sequence $\langle F_\alpha : \alpha < \omega_1\rangle$ from Fact \ref{disjoint union smooth equal reals union omega1}. Note that $\emptyset$ can serve as the $\infty$-Borel code for this sequence. Let $X \in \degrees$. $\reals^{L[X]}$ is countable. Hence there is some $\alpha < \omega_1$ so that for all $\beta > \alpha$, $\reals^{L[X]} \subseteq \reals \setminus \mathrm{WO}_\beta$. Hence for all $\beta > \alpha$, every real of $L[X]$ belongs to the single ordinal definable uncountable class of $E_\beta$. So for all $\beta > \alpha$, $X$ can not serve as the witness to the Case II assumption. This shows why the natural attempt to inject $\reals \times \omega_1$ into $\bigsqcup_{\alpha < \omega_1} \reals \slash F_\alpha \approx \reals \sqcup \omega_1$ must fail. 

However, when $\langle E_\alpha : \alpha < \omega_1\rangle$ is a sequence of smooth equivalence relations with all classes countable, this natural attempt does succeed. 

For the next theorem, one will partition $\degrees$ into two disjoint sets $\degrees = H_0^\alpha \cup H_1^\alpha$ for various $\alpha$'s. One of the two sets belongs to the ultrafilter $\mu$. The main task is to show that $H_0^\alpha$ is the one that belongs to $\mu$. Supposing for the sake of contradiction that $H_1^\alpha \in \mu$. The main trick is to code information from all the local models $\HOD_S^{L[S,X]}$ into one single ordinal by using the wellfoundedness of the ultrapower $\prod_{X \in \degrees} \omega_1\slash \mu$. This magical ordinal contains so much information that it can then be used to give an $\OD_S$ definition for the $E_\alpha$-class, $[a]_{E_\alpha}$, where $a \notin \OD_S$. Since $[a]_{E_\alpha}$ is a countable $\OD_S$ sets, this implies $a \in \OD_S$ by Fact \ref{wo od set inside hod}. This yields the desired contradition. The details follow below:

\Begin{theorem}{countable smooth R times kappa injection}
Assume $\mathsf{ZF + AD^+ + V = L(\mathscr{P}(\reals))}$. Let $\kappa \in \ON$ and $\langle E_\alpha : \alpha < \kappa\rangle$ be a sequence of smooth equivalence relations on $\reals$ with all classes countable. Then $\reals \times \kappa$ injects into $\bigsqcup_{\alpha < \kappa} \reals \slash E_\alpha$. 
\end{theorem}

\begin{proof}
By Fact \ref{natural model AD relation ordinals have infinity borel codes}, let $(S,\varphi)$ be an $\infty$-Borel code for the sequence $\langle E_\alpha : \alpha < \kappa\rangle$. By the description above, it suffices to find a single $X \in \degrees$ so that all $\alpha < \kappa$, there is some $a \in \reals^{L[S,X]}$ so that $a$ does not belong to any $E_\alpha$-component in $\bbO_S^{L[S,X]}$. 

Fix some $a \in \reals$ which is not $\OD_S$. 

Now for any $\alpha < \kappa$, let $H^\alpha_0$ be the set of $X \in \degrees$ so that $a$ does not belongs to any $E_\alpha$-component of $\bbO_S^{L[S,X]}$. Let $H^\alpha_1$ be the set of $X \in \degrees$ so that $a$ does belong to some $E_\alpha$-component of $\bbO_S^{L[S,X]}$. Since $\mu$ is an ultrafilter, either $H^\alpha_0 \in \mu$ or $H^\alpha_1 \in \mu$. Suppose that $H^\alpha_1 \in \mu$. Define $f : \degrees \rightarrow \omega_1$ by $f(X)$ is the least $\beta$ so that $a$ belongs to the $\beta^\text{th}$ $E_\alpha$-component in $\bbO_S^{L[S,X]}$ if $X \in H_1^\alpha$, and $f(X)$ is some default value otherwise. Let $F = [f]_\sim$ where $\sim$ is the $\mu$-almost equal relation. Then $[a]_{E_\alpha}$ is ordinal definable using $S$ and $F$ as parameters since $b \in [a]_{E_\alpha}$ if and only if for all representative $f$ for $F$, for a cone of $X \in \degrees$, $b$ is $E_\alpha$ related to some element in the $f(X)^\text{th}$ $E_\alpha$-component in $\bbO_S^{L[S,X]}$. By Fact \ref{wellfoundedness martin ultraproduct}, $\prod_{X \in \degrees} \omega_1 \slash \mu$ is wellfounded. Hence $F$ is essentially an ordinal. This shows that $[a]_{E_\alpha}$ is $\OD_S$. Then $[a]_{E_\alpha}$ is a countable $\OD_S$ set of reals. Fact \ref{wo od set inside hod} implies that $[a]_{E_\alpha}$ consists only of $\OD_S$ elements. But $a \notin \OD_S$ by the initial choice of $a$.. Contradiction. This shows that $H^\alpha_0 \in \mu$. 

Since for each $\alpha < \kappa$, $[a]_{E_\alpha}$ is $\OD^V_{\langle S,a \rangle}$, $[a]_{E_\alpha} \subseteq \HOD^V_{\langle S,a \rangle}$. Using the canonical wellordering of $\HOD_{\langle S,a \rangle}^V$, let $\langle b^\alpha_\epsilon : \epsilon < \delta_\alpha\rangle$ be an injective enumeration of $[a]_{E_\alpha}$, where $\delta_\alpha \in \ON$. For each $\alpha < \kappa$, $\delta_\alpha$ is less than $(2^{\aleph_0})^{\HOD_{\langle S,a\rangle}^V}$. Since $V \models \AD$, $(2^{\aleph_0})^{\HOD_{\langle S,a\rangle}^V}$ is countable in $V$. In $V$, choose a bijection of $\omega$ with $(2^{\aleph_0})^{\HOD^V_{\langle S,a\rangle}}$. This bijection then induces canonically bijections $\Sigma_\alpha : \omega \rightarrow \delta_\alpha$ for all $\alpha < \kappa$. Define $r_\alpha \in \reals$ by $r_\alpha(\langle n,k\rangle) = b^\alpha_{\Sigma_\alpha(n)}(k)$. Thus $[a]_{E_\alpha} = \{(r_\alpha)_n : n \in \omega\}$ where $(r_\alpha)_n(k) = r_\alpha(\langle n,k\rangle)$. But $\langle r_\alpha : \alpha < \kappa\rangle$ is a wellordered sequence of reals. Under $\AD$, there are only countably  many distinct $r_\alpha$'s. Thus $\{[a]_{E_\alpha} : \alpha < \kappa\}$ is a countable set.

Note that if $[a]_{E_\alpha} = [a]_{E_\beta}$, then $H^\alpha_0 = H^\beta_0$. Therefore, by countable additivity of $\mu$, $\bigcap_{\alpha < \kappa} H^\alpha_0 \in \mu$. Let $X \in \bigcap_{\alpha < \kappa} H^\alpha_0$. This is the desired degree $X$ that witnesses the Case II assumption for all $E_\alpha$. By the remarks above, this allows for the construction of an injection of $\reals \times \kappa$ into $\bigsqcup_{\alpha < \kappa} \reals \slash E_\alpha$ which completes the proof.
\end{proof}

\section{Upper Bound on Cardinality}\label{upper bound on cardinality}

\Begin{definition}{equiv relation separating family}
Let $E$ be an equivalence relation on $\reals$. Let $\mathcal{S}$ be a collection of nonempty subsets of $\reals$. $\mathcal{S}$ is a separating family for $E$ if and only if for all $x,y \in \reals$, $x \ E \ y$ if and only if for all $A \in \mathcal{S}$, $x \in A \Leftrightarrow y \in A$. 
\end{definition}

\Begin{definition}{E0 equivalence relation}
$E_0$ is the equivalence relation $\cantorspace$ defined by $x \ E_0 \ y$ if and only if $(\exists m)(\forall n \geq m)(x(n) = y(n))$. 
\end{definition}

The following is Hjorth's $E_0$-dichotomy in $\AD^+$ which generalized the classical $E_0$-dichotomy of Harrington-Kechris-Louveau \cite{Glimm-Effros-Dichotomy-for-Borel-Equivalence-Relations}.

\Begin{fact}{generalized E0 dichotomy}
(\cite{Dichotomy-for-Definable-Universe} Theorem 2.5) Assume $\ZF + \AD^+$. Let $E$ be an equivalence relation on $\reals$. Then either

(i) There is a wellordered separating family for $E$. 

(ii) There is a $\Phi : \reals \rightarrow \reals$ with the property that $x \ E_0 \ y$ if and only if $\Phi(x) \ E \ \Phi(y)$. 
\end{fact}

\begin{proof}
Note that option (ii) implies that $\reals \slash E_0$ injects in $\reals \slash E$. Suppose option (i) holds. Let $\mathcal{S} = \langle B_\alpha : \alpha < \delta\rangle$ where $\delta$ is some ordinal be the given separating family. For each $x \in \reals$, let $\Psi(x) = \{\alpha : x \in B_\alpha\}$. $\Psi$ induces an injection of $\reals \slash E$ into $\mathscr{P}(\delta)$. 

As usual in dichtomy results, there are two cases. One case yields the wellordered separating family and the other case yields an embedding. For the purpose of this paper, one is more concerned with producing the wellordered separating family. Moreover, one needs to observe that the wellordered separating family and its wellordering is produced uniformly from the $\infty$-Borel code for $E$. The following will give the argument to produce a wellordered separating family. The embedding case will be omitted as it is not relevant for this paper.

Let $(S,\varphi)$ be an $\infty$-Borel code for $E$. Regardless of the universe in consideration, $E$ will always be considered as the set defined by the $\infty$-Borel code $(S,\varphi)$. 

(Case I) For all $X \in \degrees$, for all $a,b \in \reals^{L[S,X]}$, if $\neg(a \ E \ b)$, then there is an $\OD_S^{L[S,X]}$ $C \in \bbO_S^{L[S,X]}$ which is $E$-invariant in $L[S,X]$ and $a \in C$ and $b \notin C$. 

For each $F \in \prod_{X \in \degrees} \omega_1 \slash \mu$, define $A_F$ as follows: Let $f : \degrees \rightarrow \omega_1$ be such that $f \in F$, that is $f$ is a representative of $F$. For $a \in \reals$, $a \in A_F$ if and only if on a Turing cone of $X \in \degrees$, $a$ belongs to the $f(X)^\text{th}$ $E$-invariant set in $\bbO_S^{L[S,X]}$. (If there is no $f(X)^\text{th}$ $E$-invariant $\OD_S^{L[S,X]}$-set, then let this set be $\emptyset$.) Note that $A_F$ is well defined.

$A_F$ is $E$-invariant: Suppose $a,b \in \reals$, $a \ E \ b$ and $a \in A_F$. Pick $f \in F$. Since $a \in A_F$, there is some $Z \geq_T [a \oplus b]_T$ so that for all $X \geq_T Z$, $a$ belongs to the $f(X)^\text{th}$ $E$-invariant set in $\bbO_S^{L[S,X]}$. Note $b \in L[S,X]$. $V \models a \ E \ b$ means that $V \models L[S,a,b] \models \varphi(S,a,b)$. Hence $L[S,X] \models L[S,a,b] \models \varphi(S,a,b,)$. Thus $b$ also belongs to the $f(X)^\text{th}$ $E$-invariant set in $\bbO_S^{L[S,X]}$. Hence $b \in A_F$. 

By Fact \ref{wellfoundedness martin ultraproduct}, $\prod_{X \in \degrees} \omega_1 \slash \mu$ is wellfounded. Hence $\mathcal{S} = \langle A_F : F \in \prod_{X \in \degrees} \omega_1 \slash \mu\rangle$ is a wellordered set of $E$-invariant subsets of $\reals$. 

$\mathcal{S}$ is a separating family for $E$: Suppose $a,b \in \reals$ and $\neg(a \ E \ b)$. Define $f : \degrees \rightarrow \omega_1$ as follows. Let $Z = [a \oplus b]_T$. If $X \geq_T Z$, then let $f(X)$ be the least ordinal $\alpha < \omega_1$ so that the $\alpha^\text{th}$ $E$-invariant set in $\bbO_S^{L[S,X]}$ contains $a$ but not $b$. Such an $\alpha$ exists from the Case I assumption. If $X$ is not Turing above $Z$, then let $f(X) = \emptyset$. Let $F = [f]_\sim$, where $\sim$ is the $\mu$-almost equal relation. Then $A_F \in \mathcal{S}$, $a \in A_F$, and $b \notin A_F$.

In conclusion, one has shown that $\mathcal{S}$ is a wellordered separating family for $E$.

(Case II) There is some $X \in \degrees$ and some $a,b \in \reals^{L[S,x]}$ with $\neg(a \ E \ b)$ such that there are no $E$-invariant sets $C \in \bbO_S^{L[S,X]}$ so that $a \in C$ and $b \notin C$. 

The idea is that this case assumption gives a natural condition in the forcing $\bbO_S^{L[S,X]}$ for which a perfect tree of mutual $\bbO_S^{L[S,X]}$-generic over $\HOD_S^{L[S,X]}$ below this condition serves as the desired embedding. The details can be found in \cite{Dichotomy-for-Definable-Universe} and are omitted since this case is not relevant for the rest of the paper.
\end{proof}

\Begin{theorem}{upper bound disjoint union smooth}
Assume $\mathsf{ZF + AD^+ + V = L(\mathscr{P}(\reals))}$. Let $\kappa$ be an ordinal and $\langle E_\alpha : \alpha < \kappa\rangle$ be a sequence of smooth equivalence relations on $\reals$ with all classes countable. Then there is an injection of $\bigsqcup_{\alpha < \kappa} \reals \slash E_\alpha$ into $\reals \times \kappa$. 
\end{theorem}

\begin{proof}
By Fact \ref{natural model AD relation ordinals have infinity borel codes}, let $(U,\varphi)$ be an $\infty$ code for $\langle E_\alpha : \alpha < \kappa\rangle$. Uniformly from $(U,\varphi)$, one obtains $\infty$-Borel codes $(\langle U, \alpha\rangle,\varphi_\alpha)$ for each $E_\alpha$. Since each $E_\alpha$ is smooth, $\reals \slash E_\alpha \approx \reals$. Since under $\mathsf{AD}$, $\reals \slash E_0$ does not inject into $\reals$, Case I from the proof of Fact \ref{generalized E0 dichotomy} must occur. The proof in Case I uniformly produces, from $(\langle U, \alpha\rangle,\varphi_\alpha)$, a separating family $\mathcal{S}_\alpha = \langle A_\gamma^\alpha : \gamma < \delta\rangle$ for $E_\alpha$, where $\delta$ is the ordertype of $\prod_{X \in \degrees} \omega_1 \slash \mu$. 

Let $\Phi_\alpha : \reals \rightarrow \mathscr{P}(\delta)$ be defined by $\Phi_\alpha(x) = \{\gamma : x \in A_\gamma^\alpha\}$. Since $\mathcal{S}_\alpha$ is a separating family for $E_\alpha$, $\Phi_\alpha$ is $E_\alpha$-invariant. Thus $\Phi_\alpha$ induces an injection $\tilde \Phi_\alpha$ of $\reals \slash E_\alpha$ into $\mathscr{P}(\delta)$. Define a relation $R \subseteq \kappa \times \mathscr{P}(\delta) \times \reals$ by 
$$R(\alpha,B,x) \Leftrightarrow (B = \Phi_\alpha(x)) \vee ((\forall y)(B \neq \Phi_\alpha(y)) \wedge x = \bar{0})$$
where $\bar{0}$ is the constant $0$ sequence. For each $(\alpha,B) \in \kappa \times \mathscr{P}(\delta)$, the section $R_{(\alpha,B)}$ is countable. Fact \ref{wellorderable section uniformization} implies that there is a uniformization function $F : (\kappa \times \mathscr{P}(\delta)) \rightarrow \reals$. 

Define $\Psi : \bigsqcup_{\alpha < \kappa} \reals \slash E_\alpha \rightarrow \reals \times \kappa$ by $\Psi([x]_{E_\alpha}) = (F(\alpha, \tilde\Phi_\alpha([x]_{E_\alpha})), \alpha)$. $\Psi$ is an injection.
\end{proof}

\Begin{definition}{jonsson property}
Let $X$ be a set. For $n \in \omega$, let $[X]^n_= = \{f \in {}^nX : (\forall i,j \in n)(i \neq j \Rightarrow f(i) \neq f(j)\}$. Let $[X]^{<\omega}_= = \bigcup_{n \in \omega} [X]_=^n$. 

A set $X$ has the J\'onsson property if and only if for all $f : [X]_=^{<\omega} \rightarrow X$, there is some $Y \subseteq X$ with $Y \approx X$ so that $f[[Y]_=^{<\omega}] \neq X$. 
\end{definition}

\Begin{fact}{some sets with and without jonsson}
Assume $\ZF + \AD$.

(\cite{Partition-Properties-for-Non-Ordinal-Sets}, Holshouser and Jackson) $\reals$, $\reals \sqcup \omega_1$, and $\reals \times \kappa$ where $\kappa < \Theta$ have the J\'onsson property.

(\cite{Definable-Combinatorics-Some-Borel-Equivalence-Relations-arxiv}) $\reals \slash E_0$ does not have the J\'onsson property.

(\cite{Ordinal-Definability-and-Combinatorics-Equivalence-Relations-arxiv} Fact 4.23) For any $\kappa < \Theta$, $(\reals \slash E_0) \times \kappa$ does not have the J\'onsson property.
\end{fact}

\Begin{fact}{disjoint union countable equiv pseudo-jonsson}
(\cite{Ordinal-Definability-and-Combinatorics-Equivalence-Relations-arxiv} Fact 4.13) Assume $\ZF + \AD$. Let $\kappa \in \ON$. Let $\langle E_\alpha : \alpha < \kappa\rangle$ be a sequence of equivalence relations on $\reals$ with all classes countable. Let $f : [\bigsqcup_{\alpha < \kappa}\reals \slash E_\alpha]_=^{<\omega} \rightarrow \bigsqcup_{\alpha < \kappa} \reals \slash E_\alpha$. Then there is some perfect tree $p$ so that $f[[\bigsqcup_{\alpha < \kappa} [p]\slash E_\alpha]^{<\omega}_=] \neq \bigsqcup_{\alpha < \kappa} \reals \slash E_\alpha$. 

(\cite{Ordinal-Definability-and-Combinatorics-Equivalence-Relations-arxiv} Theorem 4.15) $\reals \times \kappa$ has the J\'onsson property for all $\kappa \in \ON$. 
\end{fact}

If for every perfect tree $p$, $\bigsqcup_{\alpha < \kappa} \reals \slash E_\alpha \approx \bigsqcup_{\alpha < \kappa} [p] \slash E_\alpha$, then Fact \ref{disjoint union countable equiv pseudo-jonsson} would imply that $\bigsqcup_{\alpha < \kappa} \reals \slash E_\alpha$ has the J\'onsson property.  However, in general these two sets can not be in bijection since $\reals \slash E_0$ does not have the J\'onsson property. In this particular case, the $p$ satisfying fact \ref{disjoint union countable equiv pseudo-jonsson} is not an $E_0$-trees (see \cite{Definable-Combinatorics-Some-Borel-Equivalence-Relations-arxiv} Definition 5.2), i.e. a perfect tree with certain symmetry conditions.

Combining Theorem \ref{countable smooth R times kappa injection} and \ref{upper bound disjoint union smooth}, one can determine the cardinality of disjoint unions of quotients of smooth equivalence relations with all classes countable and show that they have the J\'onsson property.

\Begin{theorem}{cardinality disjoint union smooth}
Assume $\mathsf{ZF + AD^+ + V = L(\mathscr{P}(\reals))}$. Let $\kappa \in \ON$ and $\langle E_\alpha : \alpha < \kappa\rangle$ be a sequence of smooth equivalence relations on $\reals$ with all classes countable. Then $\bigsqcup_{\alpha < \kappa} \reals \slash E_\alpha \approx \reals \times \kappa$ and hence $\bigsqcup_{\alpha < \kappa} \reals \slash E_\alpha$ has the J\'onsson property.
\end{theorem}

\section{Disjoint Union of Quotients of Smooth Equivalence with Uncountable Classes}\label{disjoint union of quotients of smooth equivalence uncountable classes}

This section will show that in natural models of $\AD^+$ many subsets below $[\omega_1]^{<\omega_1}$ are in bijection with disjoint unions of quotients of smooth equivalence relations on $\reals$. It will be shown that any subset of $[\omega_1]^{<\omega_1}$ that contains $\reals \sqcup \omega_1$ can be written in this way. The argument is similar to the example $\langle F_\alpha : \alpha < \omega_1\rangle$ produced in the proof of Fact \ref{disjoint union smooth equal reals union omega1}. Note that each $F_\alpha$ has one uncountable equivalence class that holds the reals that are not used for coding. In the following argument, the existence of a copy of $\omega_1$ is again used to handle these classes.

Recall the distinction between smooth and weakly smooth from Definition \ref{smooth equivalence relation}. The first result of the next theorem is proved in just $\ZF + \AD$. The quotients of the weakly smooth but not smooth $E_\alpha$'s are (non-uniformly) in bijection with a countable ordinal. The uniformity of $\infty$-Borel code will be important in the argument to absorb these quotients of the weakly smooth but not smooth equivalence relations into $\omega_1$. Thus it is unclear if the second statement of the next theorem is provable in just $\AD$ or $\AD^+$.

\Begin{theorem}{characterization disjoint union smooth}
Assume $\ZF + \AD$. Suppose $X \subseteq [\omega_1]^{<\omega_1}$ and $\omega_1$ injects into $X$. Then there exists a sequence $\langle E_\alpha : \alpha < \omega_1\rangle$ of weakly smooth equivalence relations on $\reals$ so that $X$ is in bijection with $\bigsqcup_{\alpha < \omega_1} \reals \slash E_\alpha$. 

Assume $\mathsf{ZF + AD^+ + V = L(\mathscr{P}(\reals))}$. Suppose $X \subseteq [\omega_1]^{<\omega_1}$ and $\reals \sqcup \omega_1$ injects into $X$. Then there exists a sequence $\langle E_\alpha : \alpha < \omega_1\rangle$ of smooth equivalence relations on $\reals$ so that $X$ is in bijection with $\bigsqcup_{\alpha < \omega_1} \reals \slash E_\alpha$. 

Therefore, $X \subseteq [\omega_1]^{<\omega_1}$ has a sequence $\langle E_\alpha : \alpha < \omega_1\rangle$ of smooth equivalence relations such that $X \approx \bigsqcup_{\alpha < \omega_1} \reals \slash E_\alpha$ if and only if $\reals \sqcup \omega_1$ injects into $X$.
\end{theorem}

\begin{proof}
Since $\omega_1$ injects into $X$, let $X = X_0 \sqcup X_1$ where $X_1$ is in bijection with $\omega_1$. Henceforth, say $X = X_0 \sqcup \omega_1$. 

Let $\mathrm{WO}$ denote the set of reals coding wellorderings on $\omega$. Note that ${}^\omega\reals$ is in bijection with $\reals$. For $\alpha < \omega_1$, let $\mathrm{WO}^{\alpha}$ be the set of $j \in {}^\omega\mathrm{WO}$ so that for all $m \neq n$, $\mathrm{ot}(j(m)) \neq \mathrm{ot}(j(n))$ and $\sup \{\mathrm{ot}(j(n)) : n \in \omega\} = \alpha$.

For $\alpha < \omega_1$, let $X_0^{\alpha} = \{f \in X : \sup(f) = \alpha\}$. For each $j \in \mathrm{WO}^\alpha$, let $\Psi(j) \in [\alpha + 1]^{\leq\alpha}$ be the increasing enumeration of $\{\mathrm{ot}(j(n)) : n \in \omega\}$. Let $Y^{\alpha} = \Psi^{-1}[X_0^{\alpha}]$.

Define $E_{\alpha}$ on ${}^\omega \reals$ by
$$x \ E_\alpha \ y \Leftrightarrow (x \notin Y^{\alpha} \wedge y \notin Y^{\alpha}) \vee (x \in Y^{\alpha} \wedge y \in Y^{\alpha} \wedge \Psi(x) = \Psi(y))$$ 

Note that each ${}^\omega \reals \slash E_\alpha$ has one distinguished equivalence class corresponding to ${}^\omega\reals \setminus Y^{\alpha}$. Denote this class by $\star_{\alpha}$. $({}^\omega \reals \slash E_\alpha) \setminus \{\star_{\alpha}\}$ is in bijection with $X_0^{\alpha}$ in a canonical way. Thus canonically there is a bijection of $\bigsqcup_{\alpha < \omega_1} {}^\omega\reals \slash E_\alpha$ with $X_0 \sqcup \omega_1$. Also since ${}^\omega \reals \slash E_\alpha$ injects into $[\alpha + 1]^{\leq \alpha}$, which is in bijection with $\reals$, ${}^\omega\reals \slash E_\alpha$ is either in bijection with $\reals$ or is countable. (Note that this shows under $\mathsf{ZF}$ that any $X \subseteq [\omega_1]^{<\omega_1}$ which contains a copy of $\omega_1$ is a disjoint union $\bigsqcup \reals \slash E_\alpha$ where each $\reals \slash E_\alpha$ is either countable or in bijection with $\reals$.)

By Fact \ref{natural model AD relation ordinals have infinity borel codes}, there is an $\infty$-Borel code $(S,\varphi)$ for $\langle E_\alpha : \alpha < \omega_1\rangle$ in the sense of Definition \ref{infinity borel code for relation on ordinals}. As before, one can thus obtain uniformly the $\infty$-Borel code (in the ordinary sense) for each $E_\alpha$. 

Let $A = \{\alpha \in \omega_1 : |{}^\omega\reals \slash E_\alpha| = \aleph_0\}$. The argument in Case I of Theorem \ref{perfect set dichotomy} shows that uniformly in the $\infty$-Borel code for $E_\alpha$ for $\alpha \in A$, there is a wellordering of ${}^\omega\reals \slash E_\alpha$. 

Let $B = \omega_1 \setminus A$. $B \neq \emptyset$ since otherwise using the uniform wellordering of ${}^\omega \reals \slash E_\alpha$ for all $\alpha \in A = \omega_1$, one could produce a bijection of $\bigsqcup_{\alpha < \omega_1} {}^\omega\reals \slash E_\alpha$ with $\omega_1$. It was shown above that $\bigsqcup_{\alpha < \omega_1} {}^\omega\reals \slash E_\alpha$ is in bijection with $X_0 \sqcup \omega_1 = X$. However $X$ contains a copy of $\reals$. Contradiction.

Using the uniform wellordering of ${}^\omega\reals \slash E_\alpha$ for all $\alpha \in A$, the set $K = \{\star_\alpha : \alpha \in \omega_1\} \cup \bigcup_{\alpha \in A} {}^\omega\reals \slash E_\alpha$ (these two sets are not disjoint) is in bijection with $\omega_1$. If $B$ is countable, then $\bigsqcup_{\alpha < \omega_1} {}^\omega\reals \slash E_\alpha = K \cup \bigsqcup_{\alpha \in B} {}^\omega\reals \slash E_\alpha$ is in bijection with $\reals \sqcup \omega_1$. By Fact \ref{disjoint union smooth equal reals union omega1}, this set is a disjoint union of quotients of smooth equivalence relations on $\reals$. Now suppose $B$ is uncountable. Since $K \approx \omega_1$, pick a bijection of $K$ with $\{\star_\alpha : \alpha \in B\}$. Since $\bigsqcup_{\alpha < \omega_1} {}^\omega\reals \slash E_\alpha = K \sqcup \bigsqcup_{\alpha \in B} ({}^\omega\reals \slash E_\alpha) \setminus \{\star_\alpha\}$, the map that sends $K$ to $\{\star_\alpha : \alpha < \omega_1\}$ via the fixed bijection above and the identity on $\bigsqcup_{\alpha \in B} ({}^\omega\reals \slash E_\alpha) \setminus \{\star_\alpha\}$ is a bijection of $\bigsqcup_{\alpha < \omega_1} {}^\omega\reals \slash E_\alpha$ with $\bigsqcup_{\alpha \in B} {}^\omega\reals \slash E_\alpha$. 

It has been shown that $X$ is a disjoint union of quotients of smooth equivalence relations.
\end{proof}

Of course,  $\omega_1$ and $\reals$ cannot be written as a wellordered disjoint union of quotients of smooth equivalence relations. There are other examples. 

\Begin{definition}{S1 set}
Let $S_1 = \{f \in [\omega_1]^{<\omega_1} : \sup(f) = \omega_1^{L[f]}\}$. 
\end{definition}

The next result is a consequence of very powerful dichotomy results proved in \cite{The-Cardinals-Below-CountableSubsetOmega}. The following is an explicit proof of this result.

\Begin{fact}{comparison of S1}
(\cite{The-Cardinals-Below-CountableSubsetOmega}) Assume $\ZF + \AD$. $\reals$ injects into $S_1$, and $\omega_1$ does not inject into $S_1$. 
\end{fact}

\begin{proof}
For each $r \in \reals$, consider it as a subset of $\omega$. Let $\Psi(r)$ be the subset $\omega_1^{L[r]}$ consisting of $r$ and all infinite ordinals less than $\omega_1^{L[r]}$. $\Psi(r) \in S_1$ and $\Psi$ is injective as a function from $\reals$ into $S_1$. 

Suppose $\Phi : \omega_1 \rightarrow S_1$ is an injection. $\Phi$ can be coded as a subset of $\omega_1$. Note that $\sup(\{\Phi(\alpha) : \alpha \in \omega_1\}) = \omega_1$ because otherwise $\Phi$ would be injecting into $[\alpha]^{<\alpha}$ for some $\alpha < \omega_1$. The latter is in bijection with $\reals$. This would imply that there is an uncountable wellordered sequence of reals. Also note that since $L[\Phi] \models \mathsf{ZFC}$, $\omega_1^{L[\Phi]} < \omega_1$. Therefore choose some $\alpha$ so that $\sup(\Phi(\alpha)) > \omega_1^{L[\Phi]}$. However $\Phi(\alpha) \in L[\Phi]$ so $\omega_1^{L[\Phi(\alpha)]} \leq \omega_1^{L[\Phi]} < \sup(\Phi(\alpha))$. This implies that $\Phi(\alpha) \notin S_1$. Contradiction.
\end{proof}

So in $\mathsf{ZF + AD^+ + V = L(\mathscr{P}(\reals))}$, $S_1$ is not an $\omega_1$-length disjoint union of quotients of smooth equivalence relations, but $S_1 \sqcup \omega_1$ is.

\section{Almost Full Countable Section Uniformization for $[\omega_1]^\omega \times \reals$}\label{almost full countable section uniformization}

This section will show that $|[\omega_1]^\omega|$ is not below $|\bigsqcup_{\alpha < \kappa} \reals \slash E_\alpha|$ if $\langle E_\alpha : \alpha < \kappa\rangle$ is a sequence of equivalence relations on $\reals$ with all classes countable under just $\AD$. Note that by Theorem \ref{cardinality disjoint union smooth}, $\mathsf{ZF + AD^+ + V = L(\mathscr{P}(\reals))}$ is capable of proving that such a disjoint union is in bijection with $\reals \times \kappa$. It is much more evident that $[\omega_1]^\omega$ does not inject into $\reals \times \kappa$. 

\Begin{fact}{uniformization and noninjection countable sequence disjoint union}
(\cite{Ordinal-Definability-and-Combinatorics-Equivalence-Relations-arxiv}) $(\mathsf{ZF + AD})$ Let $\kappa$ be an ordinal. Let $\langle E_\alpha : \alpha < \kappa\rangle$ be a sequence of equivalence relations on $\reals$. Let $\Phi : [\omega_1]^{\omega} \rightarrow \bigsqcup_{\alpha < \kappa} \reals \slash E_\alpha$. Let $R \subseteq [\omega_1]^\omega \times \reals$ be defined by $R(f,x) \Leftrightarrow x \in \Phi(f)$. If there is a $Z \subseteq [\omega_1]^\omega$ with $Z \approx [\omega_1]^\omega$ and a $\Lambda : Z \rightarrow \reals$ so that for all $f \in Z$, $R(f,\Lambda(f))$, then $\Phi$ is not an injection.
\end{fact}

\begin{proof}
See \cite{Ordinal-Definability-and-Combinatorics-Equivalence-Relations-arxiv} Fact 4.19 and the subsequent discussions.
\end{proof}

Using some of the ideas above, one can prove in $\AD^+$, the (full) countable section uniformization for relations on $\reals \times [\omega_1]^\omega$: For every $R \subseteq [\omega_1]^\omega \times \reals$ such that $R_f = \{x \in \reals : R(f,x)\}$ is nonempty and countable for all $f \in [\omega_1]^\omega$, there is a uniformization function for $R$. Then Fact \ref{uniformization and noninjection countable sequence disjoint union} gives the following result:

\Begin{fact}{AD+ no injection countable sequence into disjoint union}
(\cite{Ordinal-Definability-and-Combinatorics-Equivalence-Relations-arxiv} Fact 4.20) Assume $\ZF + \AD^+$. Let $\kappa$ be an ordinal and $\langle E_\alpha : \alpha < \kappa\rangle$ be a sequence of equivalence relations on $\reals$ such that each $E_\alpha$ has all classes countable. Then there is no injection $\Phi : [\omega_1]^\omega \rightarrow \bigsqcup_{\alpha < \kappa} \reals \slash E_\alpha$. 
\end{fact}

For the J\'onsson property, often uniformization on a sufficiently big set is enough for the desired result. $\AD$ can prove an almost full uniformization result for relation on $\reals \times \reals$. (For example, $\AD$ proves comeager uniformization.) Fact \ref{uniformization and noninjection countable sequence disjoint union} only requires that one can uniformize relations on $[\omega_1]^\omega \times \reals$ on a set $Z \subseteq [\omega_1]^\omega$ that has the same cardinality as $[\omega_1]^\omega$. In general, this is impossible. By Theorem \ref{characterization disjoint union smooth}, every cardinal below $[\omega_1]^{<\omega_1}$ that contains a copy of $\reals \sqcup \omega_1$, for example $[\omega_1]^\omega$, can be written as a disjoint union of smooth equivalence relations. If this almost full uniformization exists for all relations on $[\omega_1]^\omega \times \reals$, then Fact \ref{uniformization and noninjection countable sequence disjoint union} would imply that there is no injection of $[\omega_1]^\omega$ into $[\omega_1]^\omega$. 

The following will show in $\AD$ alone that one can prove almost full uniformization for relations on $[\omega_1]^\omega \times \reals$ with all sections countable. This will suffice to proves the statement of Fact \ref{AD+ no injection countable sequence into disjoint union} in $\AD$ alone.

\Begin{definition}{coding sequence}
Fix a recursive bijection $\langle \cdot, \cdot \rangle : \omega \times \omega \rightarrow \omega$. For $x \in \reals$, let $(x)_m \in \reals$ be defined by $(x)_m(n) = x(\langle m,n\rangle)$. Using the pairing function, one can also code relations on $\omega$ of various arity as a subset of $\omega$. 

Let $\mathrm{WO}$ denote the reals coding wellorderings on $\omega$. 

Let $\mathrm{WO}^\omega$ be the set of reals $x$ so that for all $n \in \omega$, $(x)_n \in \mathrm{WO}$ and for all $m < n$, $\mathrm{ot}((x)_m) < \mathrm{ot}((x)_n)$. If $x \in \mathrm{WO}^\omega$, then let $f_x \in [\omega_1]^\omega$ be defined by $f_x(n) = \mathrm{ot}((x)_n)$. (Every element of $[\omega_1]^{\omega}$ has a code in $\mathrm{WO}^\omega$). 

Fix throughout,  $W \in \mathrm{WO}$ to be a recursive wellordering of ordertype $\omega \cdot \omega$. In context, $\alpha < \omega \cdot \omega$ will refer to the element of $\omega$ which corresponds to the ordinal $\alpha$ according to the wellordering $W$. Similarly in context, $<$ will refer to the wellordering given by $W$.

Let $\mathrm{WO}^{\omega\cdot\omega}$ denote the set of $x \in \reals$ so that for all $n \in \omega$, $(x)_n \in \mathrm{WO}$, and for all $\alpha,\beta \in \omega\cdot\omega$, $\alpha < \beta$ if and only if $\mathrm{ot}((x)_\alpha) < \mathrm{ot}((x)_\beta)$. (Here $\alpha$ and $\beta$ refer to the natural numbers corresponding to $\alpha$ and $\beta$, respectively, according to $W$.) If $x \in \mathrm{WO}^{\omega\cdot\omega}$, let $g_x \in [\omega_1]^{\omega\cdot\omega}$ be defined by $g_x(\alpha) = \mathrm{ot}((x)_\alpha)$. Every element of $[\omega_1]^{\omega\cdot\omega}$ is of the form $g_x$ for some $x \in \mathrm{WO}^{\omega\cdot\omega}$.

If $x \in \mathrm{WO}^{\omega\cdot\omega}$, let $h_x \in [\omega_1]^\omega$ be defined by $h_x(n) = \sup\{g_x(\omega \cdot n + i) : i \in \omega\}$.

If $g \in [\omega_1]^{\omega\cdot\omega}$, then let $\tilde g \in [\omega_1]^\omega$ be defined by $\tilde g(n) = \sup\{g(\omega \cdot n + i) : i \in \omega\}$. 

If $X \subseteq \omega_1$, then $\mathrm{WO}_X$ be the $x \in \mathrm{WO}$ so that $\mathrm{ot}(x) \in X$. $\mathrm{WO}^\omega_X$ and $\mathrm{WO}^{\omega\cdot\omega}_X$ are defined as above with $\mathrm{WO}$ replaced by $\mathrm{WO}_X$. 

Suppose $C \subseteq \omega_1$ is closed and unbounded. Let 
$${}_\omega C = \{\sup\{C(\omega\cdot \alpha + i) : i \in \omega\} : \alpha < \omega_1\}.$$
A $C$-witness for a $f \in [{}_\omega C]^\omega$, is a function $g \in [C]^{\omega\cdot\omega}$ so that $\tilde g = f$. A $C$-code for a function $f \in [{}_\omega C]^\omega$ is an $x \in \mathrm{WO}^{\omega\cdot\omega}_C$ so that $h_x = f$. So a code in $\mathrm{WO}^{\omega\cdot\omega}_C$ for any $C$-witness for $f$ is a $C$-code. Every $f \in [{}_\omega C]^\omega$ has a $C$-code.
\end{definition}

\Begin{theorem}{representation of relation on seq}
$(\AD)$ Let $R \subseteq [\omega_1]^{\omega} \times \reals$ be such that for all $f \in [\omega_1]^\omega$, $R_f := \{x \in \reals : (f,x) \in R\}$ is nonempty. 

There exists a $\sigma \in \reals$ and some closed and unbounded $C \subseteq \omega_1$, so that for all $x \in \mathrm{WO}_C^{\omega\cdot\omega}$ so that $h_x \in [{}_\omega C]^\omega$, there is some $z \leq_T x \oplus \sigma$ so that $R(h_x, z)$. Moreover there is some formula $\varphi$ so that for all $x \in \mathrm{WO}^{\omega\cdot\omega}_C$ with $h_x \in [{}_\omega C]^\omega$ and $z \in \reals$, $L[\sigma,x,z] \models \varphi(\sigma,x,y)$ implies $R(h_x,z)$. 

Assume further that for all $f \in [\omega_1]^\omega$, $|R_f| \leq \aleph_0$. There exists some uncountable $X \subseteq \omega_1$ and function $\Psi$ which uniformizes $R$ on $[X]^\omega$: For $f \in [X]^\omega$, $R(f,\Psi(f))$. 
\end{theorem}

\begin{proof}
The first half of this argument is similar to Martin's proof of the partition relation $\omega \rightarrow (\omega_1)^\omega_2$. The second half is similar to the proof of Woodin's countable section uniformization for relations on $\reals \times \reals$ as exposited in \cite{Ramsey-Ultrafiler-and-Countable-to-One-Uniformation}.

Consider the following game: Player 1 and 2 take turns playing integers. Player 1 produces a real $x$. Player 2 produces two reals $y$ and $z$. 
$$\begin{matrix}
\mathrm{Player} \ 1  & x \\
\mathrm{Player} \ 2  & y,z
\end{matrix}$$

(As before, ordinals $\alpha < \omega\cdot\omega$ are considered natural numbers according to the fixed wellordering $W$ when required in context.) 

(Case A) If there is a least $\alpha < \omega\cdot\omega$ so that 

\noindent (i) $(x)_\alpha \notin \mathrm{WO}$ or $\mathrm{ot}((x)_\alpha) \leq \sup\{\mathrm{ot}((x)_{\beta}),\mathrm{ot}((y)_{\beta}) : \beta < \alpha\}$ or 

\noindent (ii) $(y)_\alpha \notin \mathrm{WO}$ or $\mathrm{ot}((y)_\alpha) \leq \sup\{\mathrm{ot}((x)_{\beta}),\mathrm{ot}((y)_{\beta}) : \beta < \alpha\}$.

\noindent Player 2 wins if and only if $(i)$ holds. 

(This means that one of the two players fails to ensure that every section of its own real codes a wellordering or fails to produce a wellordering larger than the wellordering of all the previous sections of both reals. In this case Player 2 wins if and only if Player 1 is the first to fail in this manner.)

(Case B) Suppose there is no such $\alpha$ as above. Define $h \in [\omega_1]^\omega$ by $h(\alpha) = \sup\{\mathrm{ot}((x)_\alpha),\mathrm{ot}((y)_\alpha)\}$. Player 1 wins if $\neg R(\tilde h, z)$. 

This completes the definition of the payoff set of the game.
\newline
\newline\indent Claim 1: Player 1 does not have a winning strategy in this game. 

To see this. Suppose $\sigma$ is a winning strategy for Player 1. For each $\alpha < \omega \cdot \omega$ and $\beta < \omega_1$, let $B^\alpha_\beta$ be the set of $r \in \mathrm{WO}$ so that there exists some $y \in \reals$ and $z \in \reals$ so that for all $\gamma < \alpha$, $(y)_\gamma \in \mathrm{WO}$, $\sup\{\mathrm{ot}((y)_\gamma) : \gamma < \alpha\} < \beta$, and if $x$ is the result of Player 1 in $\sigma * (y,z)$, then $r = (x)_\alpha$. $B^\alpha_\beta$ is $\analytic$ (using a code for $\beta$ as a parameter). By the boundedness principle, there is some $\delta_\beta^\alpha < \omega_1$ so that for all $r \in B_\beta^\alpha$, $\mathrm{ot}(r) < \delta_\beta^\alpha$. Let $C$ be the set of $\eta < \omega$ so that for all $\alpha < \omega \cdot \omega$ and $\beta < \eta$, $\delta_\beta^\alpha < \eta$. 

Let $f \in [{}_\omega C]^\omega$. Pick some $g \in [C]^{\omega\cdot\omega}$ so that $\tilde g = f$. Let $y \in \mathrm{WO}^{\omega\cdot\omega}_C$ be such that $g_y = g$. Let $z \in R_f$. Play $\sigma * (y,z)$. Let $x$ be the response produced by Player 1 according to $\sigma$. Define $h(\alpha) = \sup\{\mathrm{ot}((x)_\alpha), \mathrm{ot}((y)_\alpha)\}$. By definition of $C$, $\mathrm{ot}((x)_\alpha) < \mathrm{ot}((y)_\alpha)$. Thus $\tilde h = \tilde g = f$. Then $R(\tilde h, z)$. Player 2 won. This contradicts $\sigma$ being a winning strategy for Player 1. This completes the proof of Claim 1.
\newline
\newline\indent Now suppose that $\sigma$ is a winning strategy for Player 2. For each $\alpha < \omega \cdot \omega$ and $\beta < \omega_1$, let $B^\alpha_\beta$ be the set of $r \in \mathrm{WO}$ so that there exists some $x \in \reals$, so that for all $\gamma \leq \alpha$, $(x)_\gamma \in \mathrm{WO}$, $\sup\{\mathrm{ot}((x)_\gamma) : \gamma < \alpha\} < \beta$, and if $(y,z)$ is the response of Player 2 from $x * \sigma$, then $r = (y)_\alpha$. Each $B_\beta^\alpha \subseteq \mathrm{WO}$ and is $\analytic$. By the boundedness lemma, there is a least ordinal $\delta_\beta^\alpha$ so that for all $r \in B_\beta^\alpha$, $\mathrm{ot}(r) < \delta^\alpha_\beta$. Let $C$ be the closed and unbounded set of $\eta$ so that for all $\alpha < \omega\cdot\omega$ and $\beta < \eta$, $\delta^\alpha_\beta < \eta$. 

Now suppose $x \in \mathrm{WO}_C^{\omega\times\omega}$ be such that $h_x \in [{}_\omega C]^\omega$. Use the player 2 stategy $\sigma$ to play against $x$ to produce the play $x * \sigma$. Let $(y,z)$ be Player 2's response from the play $x * \sigma$. Let $h(\alpha) = \sup\{\mathrm{ot}((x)_\alpha),\mathrm{ot}((y)_\alpha)\}$. Using the definition of $C$ as before, $\tilde h = h_x$. Since $\sigma$ is winning for Player 2, one has that $R(h_x,z)$. Note that $z \leq_T x \oplus \sigma$. From this description, one can allow $\varphi(\sigma,x,z)$ to be the formula that asserts that there is some $y$ so that $(y,z)$ is Player 2 response in the play $x * \sigma$.

Now to prove the uniformization result: Assume that for all $f \in [\omega_1]^\omega$, $|R_f| \leq \aleph_0$. Let $C$ be the club set from above. By a result of Solovay (\cite{Infinitary-Combinatorics-and-the-Axiom-of-Determinateness} Lemma 2.8), there is some $w \in \reals$ so that $C$ is definable in $L[w]$ from some fixed formula using $w$. Hence ${}_\omega C$ is also definable in $L[w]$. Now suppose that $f \in [{}_\omega C]^\omega$. Let $x \in \mathrm{WO}^{\omega\cdot\omega}_C$ be such that $h_x = f$. Suppose $(y,z)$ is Player 2's response in the play $x * \sigma$.

Fix some $X \geq_T [x]_T$. In $L[\sigma,w,f,X]$, define the condition $p \in {}_2\bbO_{\sigma,w,f,X}^{L[\sigma,w,X]}$ by
$$p = \{(a,b) \in \reals^2 : L[\sigma,w,f,a,b] \models \psi(\sigma,w,a,b)\}$$
where $\psi(\sigma,w,f,a,b)$ asserts that $a$ is a $C$-code for $f$ and $\varphi(\sigma,a,b)$. Note that one uses the real $w$ to speak about $C$ in $L[\sigma,w,f,a,b]$. Note that $p \neq \emptyset$ since $(x,z) \in p$. (Note that $z \leq_T \sigma \oplus x$ and hence $z \in L[\sigma,w,f,X]$.) This shows that $p$ is indeed a condition in ${}_2\bbO_{\sigma,w,f}^{L[\sigma,w,f,X]}$. 

Let $\tau$ denote the canonical name for the generic element of $\reals^2$ added by ${}_2\bbO_{\sigma,w,f}^{L[\sigma,w,f,X]}$. 
\newline
\newline\indent Claim 2: There is a dense set of condition below $p$ that determines the value of the second coordinate of $\tau$.

To prove this: Suppose not. Let $p' \leq p$ be some condition so that no $q \leq p'$ determines the second coordinate of $\tau$. Since $\AD$ holds and $\HOD_{\sigma,w,f}^{L[\sigma,w,f,X]} \models \AC$,  $(\mathscr{P}({}_2\bbO^{L[\sigma,w,f,X]}_{\sigma,w,f}))^{\HOD_{\sigma,w,f}^{L[\sigma,w,f,X]}}$ is countable in $V$. In $V$, let $\langle D_n : n \in \omega\rangle$ enumerate all the dense open subsets of ${}_2\bbO_{\sigma,f,w}^{L[\sigma,f,w]}$ that belong to $\HOD_{\sigma,f,w}^{L[\sigma,f,w,X]}$. Let $\pi_2: \reals^2 \rightarrow \reals$ be the projection onto the second coordinate.

Let $p_\emptyset \leq p'$ be the least element below $p$ meeting $D_0$ according to the canonical wellordering of $\HOD_{\sigma,w,f}^{L[\sigma,w,f,X]}$. Let $m_\emptyset = 0$. Suppose for some $\sigma \in \finBinarySequence$, $p_\sigma$ and $m_\sigma$ have been defined. Since $p_\sigma \leq p'$, no condition extending $p_\sigma$ can determine $\pi_2(\tau)$. Thus there is some $N > m_\sigma$ and some least pair $p_0,p_1 \leq p_\sigma$ so that, $p_0,p_1 \in D_{|\sigma| + 1}$, $p_0$ and $p_1$ both decides $\pi_2(\tau) \upharpoonright N$ and $p_i \forces \pi_2(\tau)(\check N) = \check i$ (that is, decides the value at $N$ differently). Let $m_{\sigma\hat{\ }i} = N + 1$ and $p_{\sigma\hat{\ }i} = p_i$ for both $i \in 2$. This produces a sequence $\langle p_\sigma : \sigma \in \finBinarySequence\rangle$. 

For each $r \in \reals = \cantorspace$, let $G_r$ be the upward closure of $\{p_{r\upharpoonright n} : n \in \omega\}$. $G_r$ is a ${}_2\bbO_{\sigma,w,f}^{L[\sigma,w,f,X]}$-generic filter over $\HOD_{\sigma,w,f}^{L[\sigma,w,f,X]}$. Also by construction, if $r \neq s$, then $\pi_2(\tau)[G_r] \neq \pi_2(\tau)[G_s]$. For all $r \in \reals$, $p \in G_r$. Since $p$ is a condition of the form to apply Fact \ref{vopenka theorem}, one has that $\HOD_{\sigma,w,f}^{L[\sigma,w,f,X]}[G_r] \models L[\sigma,w,f,\pi_1(\tau[G_r]),\pi_2(\tau[G_r])] \models \psi(\sigma,w,f,\pi_1(\tau[G_r]),\pi_2(\tau[G_r))$. By the absoluteness of the coding, $\pi_1(\tau[G_f])$ is a $C$-code for $f$ in $V$. By the property of the formula $\varphi$ (namely its upward absoluteness), one has that $R(f,\pi_2(\tau[G_r]))$ holds in $V$. Thus it has been shown that for all $r \in \reals$, $\pi_2(\tau[G_r]) \in R_f$. This contradicts $|R_f| \leq \aleph_0$. Claim 2 has been proved.
\newline
\newline\indent Claim 3: $\HOD_{\sigma,w,f}^{L[\sigma,w,f,X]} \cap R_f \neq \emptyset$. 

To prove this: Since $X \geq_T [x]_T$, the real $x$ and its associated $z \leq_T x \oplus \sigma$ (picked above) belong to $L[\sigma,f,w,X]$. By Fact \ref{vopenka theorem}, let $G$ be the ${}_2\bbO_{\sigma,w,f}^{L[\sigma,w,f,X]}$-generic filter over $\HOD_{\sigma,w,f}^{L[\sigma,w,f,X]}$ so that $\tau[G] = (x,z)$. Note that $p \in G$. So Fact \ref{vopenka theorem} implies that $R(f,z)$. Let $D$ be the dense set below $p$ from Claim 2. By genericity, $D \cap G$. Hence there is some $q \in D \cap G$. Since $q$ completely determines $\pi_2(\tau)$, one has that for all $i \in 2$, $z(n) = i$ if and only if $q \forces \pi_2(\tau)(n) = i$. Since $q \in {}_2\bbO_{\sigma,w,f}^{L[\sigma,w,f,X]}$, $q$ is essentially an ordinal. This shows that $z$ is $\OD_{\sigma,w,f}^{L[\sigma,w,f,X]}$. Claim 3 has been proved. 
\newline
\newline\indent It has been shown that for all $f \in [{}_\omega C]^\omega$, there is a cone of $X \in \degrees$ so that $R_f \cap \HOD_{\sigma,w,f}^{L[\sigma,w,f,X]} \neq \emptyset$. 
\newline
\newline\indent
Claim 4: There is a function $\Phi : [{}_\omega C]^\omega \rightarrow \reals$ that uniformizes $R$ on $[{}_\omega C]^\omega$. 

Fix an $f \in [{}_\omega C]^\omega$. Let $Y \in \degrees$ be a base of a cone of $X \in \degrees$ so that $R_f \cap \HOD_{\sigma,w,f}^{L[\sigma,w,f,X]} \neq \emptyset$. For each $n \in \omega$ and $i \in 2$, let $E_n^i$ be the set of $X \in \degrees$ such that $X \geq_T Y$ and if $z \in \reals$ is the least element of $\HOD_{\sigma,w,f}^{L[\sigma,w,f,X]}$ belonging to $R_f$, then $z(n) = i$. Since $E_n^0 \cap E_n^1 = \emptyset$, $E_n^0 \cup E_n^1$ is the set of all degrees above $X$, and Martin's measure $\mu$ is an ultrafilter, there is some $a_n$ so that $E_n^{a_n} \in \mu$. Let $\Phi(f) \in \reals$ be defined by $\Phi(f)(n) = a_n$. Using $\mathrm{AC}_\omega^\reals$, one can find a sequence of reals $\langle x_i : i \in \omega\rangle$ so that the cone above $[x_i]_T$ is contained in $E_n^{a_n}$. If $X \geq_T [\bigoplus_{i \in \omega} x_i]_T$, then $\Phi(f)$ is the $\HOD_{\sigma,w,f}^{L[\sigma,w,f,X]}$-least element of $R_f \cap \HOD_{\sigma,w,f}^{L[\sigma,w,f,X]}$. In particular, $R(f,\Phi(f))$. 
\end{proof}

\Begin{corollary}{no injection seq into disjoint union}
$(\ZF + \AD)$ Let $\langle E_\alpha : \alpha < \kappa\rangle$ be a sequence of equivalence relations on $\reals$ with all classes countable, then $[\omega_1]^\omega$ does not inject into $\bigsqcup_{\alpha < \kappa} \reals \slash E_\alpha$. 
\end{corollary}

\begin{proof}
This follows from Fact \ref{uniformization and noninjection countable sequence disjoint union} and Theorem \ref{representation of relation on seq}.
\end{proof}

The following result states that under $\AD$, given any arbitrary function $\Phi : [\omega_1]^\omega \rightarrow \reals$, one can find two reals $\sigma$ and $w$ and a set $X \subseteq \omega_1$ with $|X| = \omega_1$ so that for all $f \in [X]^\omega$, $\Phi(f)$ is constructible from $\sigma$, $w$, and $f$.

\Begin{theorem}{constructibility value of function almost everywhere}
$(\ZF + \AD)$. Let $\Phi : [\omega_1]^\omega \rightarrow \reals$ be a function. Then there is an uncountable $X \subseteq\omega_1$, reals $\sigma,w \in \reals$, and a formula $\phi$ so that for all $f \in [X]^\omega$, $\Phi(f) \in L[\sigma,w,f]$ and for all $z \in \reals$, $z = \Phi(f)$ if and only if $L[\sigma,w,f,z] \models \phi(\sigma,w,f,z)$. 
\end{theorem}

\begin{proof}
Treating $\Phi$ as a relation, run the same argument as in Theorem \ref{representation of relation on seq}. This produces the Player 2 winning strategy $\sigma$ and a closed and unbounded set $C$. Let $X = {}_\omega C$. Using $\AD$ and a result of Solovay (\cite{Infinitary-Combinatorics-and-the-Axiom-of-Determinateness} Lemma 2.8), let $w$ be a real so that $C$ can be defined in $L[w]$. 

Let $f \in [X]^\omega$. A $C$-code for $f$ exists in any $\mathrm{Coll}(\omega,\sup(f))$-generic extension of $L[\sigma,w,f]$. Let $\tau$ be a homogeneous name for a $C$-code for $f$. Homogeneous here means that for any formula $\varsigma$, either $L[\sigma,w,f] \models 1_{\mathrm{Coll}(\omega,\sup(f))} \forces \varsigma(\tau)$ or $L[\sigma,w,f] \models 1_{\mathrm{Coll}(\omega,\sup(f))} \forces \neg\varsigma(\tau)$. (The real $w$ is needed to speak about $C$ in $L[\sigma,w,f]$.) Define a real $z$ by: $n \in z$ if and only $1_\mathrm{Coll}(\omega,\sup(f))$ forces that when Player 1 plays $\sigma$ against $\tau$, if $(y',z')$ is the response from Player 2, then $n \in z'$. By the homogeneity of $\mathrm{Coll}(\omega,\sup(f))$ and the fact that $\Phi$ is a function, one can show that $1_{\mathrm{Coll}(\omega,\sup(f))}$ either forces the statement above or its negation. By the definability of the forcing relation $z \in L[\sigma,w,f]$. By definition of the game, $z = \Phi(f)$. 

The description above also provides the formula $\phi$.
\end{proof}

\bibliographystyle{amsplain}
\bibliography{references}

\end{document}